\documentclass[12pt]{article}
\usepackage{latexsym}
\usepackage{epsfig}
\usepackage{enumerate}
\usepackage{amsmath,amsthm,amssymb}
\usepackage[active]{srcltx}
\usepackage{epstopdf}

\usepackage[usenames]{color}
\usepackage{graphicx}

\definecolor{brown}{cmyk}{0, 0.72, 1, 0.45}
\definecolor{duck-shit}{cmyk}{0.3,0,1,0.4}
\definecolor{my-cyan}{cmyk}{1,0.4,0,0.2}
\definecolor{my-green}{cmyk}{.8,0,1,0}

\definecolor{grey}{gray}{0.5}

\newcounter{rot}

\newcommand{\bignote}[1]{   }

\newcommand{\xbignote}[1]{}

\newcommand{\bM}[1]{{\bf M#1}}

\def\a{\alpha} \def\b{\beta} \def\d{\delta} 
\def\e{\epsilon} \def\f{\phi} \def\F{{\Phi}}  \def\g{\gamma}
\def\G{\Gamma}  
\def\z{\zeta} \def\th{\theta}  \def\Th{\Theta}  \def\l{\lambda}
 \def\m{\mu} \def\n{\nu} 
\def\r{\rho}  \def\s{\sigma} 
\def\t{\tau} \def\om{\omega}

\def\cT{{\cal T}}

\newtheorem{theorem}{Theorem}
\newtheorem{lemma}[theorem]{Lemma}

\newcommand{\proofstart}{{\bf Proof\hspace{2em}}}
\newcommand{\proofend}{\hspace*{\fill}\mbox{$\Box$}}
\newcommand{\W}[1]{W}

\newcommand{\ooi}{(1+o(1))}
\newcommand{\sooi}{\mbox{${\scriptstyle (1+o(1))}$}}
\newcommand{\ul}[1]{\mbox{\boldmath$#1$}}

\newcommand{\wh}[1]{\widehat{#1}}
\newcommand{\wt}[1]{\widetilde{#1}}

\newcommand{\rdup}[1]{{\left\lceil #1 \right\rceil }}
\newcommand{\rdown}[1]{{\left\lfloor #1\right \rfloor}}

\newcommand{\beq}[1]{\begin{equation}\label{#1}}
\newcommand{\eeq}{\end{equation}}
\newcommand{\blem}[1]{\begin{lemma}\label{#1}}
\newcommand{\elem}{\end{lemma}}
\newcommand{\bthm}[1]{\begin{theorem}\label{#1}}
\newcommand{\ethm}{\end{theorem}}

\newcommand{\brac}[1]{\left(#1\right)}
\newcommand{\sbrac}[1]{\left({\scriptstyle #1}\right)}

\newcommand{\bfrac}[2]{\left(\frac{#1}{#2}\right)}

\newcommand{\rai}{\rightarrow \infty}
\newcommand{\ra}{\rightarrow}

\newcommand{\set}[1]{\left\{#1\right\}}

\def\E{\mbox{{\bf E}}}
\def\Pr{\mbox{{\bf Pr}}}
\def\whp{w.h.p.}

\newcommand{\ignore}[1]{}

\begin{document}

\makeatletter
\title{The height of random $k$-trees and related branching processes}

\author{Colin Cooper\thanks{Department of  Informatics,
King's College, University of London, London WC2R 2LS, UK.
Supported in part by EPSRC grant EP/J006300/1}\and Alan
Frieze\thanks{Department of Mathematical Sciences, Carnegie Mellon
University, Pittsburgh PA15213, USA. Supported in part by NSF grant CCF0502793.}
\and Ryuhei Uehara
\thanks{School of Information Science, JAIST, Asahidai 1-1, Nomi, Ishikawa 923-1292, Japan}}

\maketitle \makeatother

\begin{abstract}
We consider the height of random $k$-trees and $k$-Apollonian networks.
These random graphs are not really trees, but instead have
a tree-like structure. The height will be the maximum distance of
a vertex from the root. We show that w.h.p. the height of random
$k$-trees and $k$-Apollonian networks is asymptotic to $c\log t$, where $t$ is the number of vertices, and
$c=c(k)$ is given as the solution to a transcendental equation. The
equations are slightly different for the two types of process. In the limit as $k \rai$
the height of both processes is asymptotic to $\log t/(k \log 2)$.
\end{abstract}

\section{Introduction}

We give a general method for obtaining the height of tree-like random processes,
and illustrate the method by application to random $k$-trees and Apollonian networks.

The processes that we consider generate a sequence of graphs
$G(t), t \ge 0$ where $G(t)$ is obtained from
$G(t-1)$ by the addition of an extra vertex in some way. The initial
structure is
a $k$-clique with a distinguished vertex
$v$, which we use as the root vertex.
Of course, $G(t)$ is not necessarily a tree but it is convenient to adopt the
terminology.

The height of a vertex $u$ in $G(t)$ is its graph distance $d(v,u)$
from the root vertex $v$.
The height $h(G(t))$ of $G(t)$ is the maximum height of one of its vertices.
By considering the breadth first search tree $T_v$ rooted at $v$ we can
partition the vertices $u \in V$
into sets $L_i, i=0,1,...,h_v$ based on the distance $i=d(v,u)$ from the root vertex $v$.
We refer to sets $L_i$
as the level sets of the BFS tree. The height of $G(t)$ is thus the height of $T_v$
in the usual sense.

The general properties of random $k$-trees  have been investigated by
several authors including \cite{CU}, \cite{Gao}, \cite{PS}.
In particular, in  an earlier study into the small world properties of random $k$-trees,
 Cooper and Uehara \cite{CU} found experimentally that the  diameter of such trees
was a rapidly decreasing function of $k$.
The main result of this paper, given in Theorem \ref{Th1}, is that the diameter $D_k(t)$
of a random $t$ vertex $k$-tree satisfies
\[
\lim_{t \rai \atop k \rai} \frac{k}{\log t}\;D_k(t) = \frac{2}{\log 2}.
\]

\paragraph{The height of branching processes: Related work.}

The work in this area is so extensive it is impossible to summarize concisely.
As our interest
lies in the area of discrete random structures we must necessarily restrict
our discussion to those
authors who have had a direct influence on us, and on the techniques we use
in this paper.
Foremost among these are the works of Broutin and Devroye \cite{BroDev},
Devroye \cite{Dev1, Dev1a, Dev2}, Kingman \cite{King} and Pittel \cite{Pit}.
The formulation in these papers  differs from the discrete context
in which a new vertex is added at each step $t$,
but the end product is the same.
The basic model is a continuous time  reproductive process  in which the reproductive
rate $\l(j)$ of the parent depends on the number of offspring $j$.
 Each child independently reproduces according to the same process.
Such processes are known as a Crump-Mode-Jagers process (see Devroye \cite{Dev2}).
The paper of Kingman \cite{King} concerns the time $B_N$ of the first
birth in the $N$-th generation of
an age dependent reproductive process of Crump-Mode type, with a proof
that  $B_N/N \ra c$ as $N \rai$. To determine the constant $c$, the work
uses the Cram\'er function of the process, which (crudely) is an optimization
of the logarithm of the moment generating function of the
distribution of reproduction waiting times. A full description of Cram\'er
functions can be found
in \cite{BroDev}.
Pittel \cite{Pit} applied Kingman's result to  a branching process in which the number of
children born to a parent within $T$ steps is negative
exponential with linear population dependent rate $\l(j)= a j +1$ for some $a\ge 0$.
This serves as  a  model of random recursive trees where a vertex chooses
its parent $v$ with probability proportional to $a d(v)+1$,
where $d(v)$ is the out-degree of $v$ in the orientation of the tree away from the root.
The cases of random and
preferential attachment trees follow from setting $a=0$, and $a=1$ respectively.
The general solution being that
the height $h_n$ of an $n$ vertex tree satisfies $h_n \sim c \log n$,
where $c=1/((a+1)\g)$ and $\g$ is the positive
root of $a\g+ \log \g +1=0$. (We use $A_n\sim B_n$ to denote $A_n=(1+o(1))B_n$ as
$n\to\infty$).
In a sequence of papers, Devroye and  Broutin and Devroye,  develop a general approach in
which the central structure  is an infinite tree with branching
factor $b$ and a pair of independent random variables $(Z,E)$ on the edges.
The random variable $Z$ measures
increments in weighted  height, and $E$ measures the delay between birth of the parent and birth of the child.
Typically $E$ would be negative exponential rate 1. The height of a vertex $u$ is the sum of the $Z$ entries
on the path from the vertex $u$ to the root $v$. If $h_T$ is the maximum height of the subtree at time $T$,
then $h_T/T \ra c$, where $c$ is the maximum along a particular curve of an identity based on the
Cram\'er functions of $Z$ and $E$, thus extending the original proof of Kingman \cite{King}.
A complete explanation of the technique, and a wide range of supporting examples are given in \cite{BroDev}.
In general, for branchings based on the minimum of exponential waiting times with varying rate parameter,
the time $T$  and population size $n(T)$ are related by $T= \Th(\log n)$. The need to obtain the explicit
constant somewhat complicates the discussion.

\paragraph{The height of random $k$-trees.}

In the area of graph algorithms, $k$-trees form a well known graph
class that generalize trees and play an important role in the study of graph
minors (see \cite{Bod93,Bod98} for further details). One definition (among many)  of
$k$-trees is the following:

For any fixed positive integer $k\ge 2$,
(i) A complete graph $K_k$ of $k$ vertices is a $k$-tree,
(ii) For a $k$-tree $G$ of $t$ vertices, a new $k$-tree $G'$ of $t+1$ vertices is obtained by
adding a new vertex $v$ incident to a clique of size $k-1$ in $G$.
When $k=2$, this process forms a tree, by extending a chosen vertex.
When $k=3 $, the process forms a tree of triangles by extending the chosen edge, and so on.

The preferential attachment method for generating  random tree processes extends the graph by attaching a new vertex
to an existing vertex
chosen with probability proportional to its degree.
This  is equivalent to choosing to attach to a random end point of a random edge, i.e.
a random $K_1$ of a random $K_2$.
The $k$-tree process described below generalizes this approach,
in that we attach the new vertex to a random $(k-1)$-clique of a random $k$-clique.

A random $k$-tree $G_k(t)$, $t \ge k$, is obtained as follows.
Start with $G_k(1)$, a $k$-clique $C_1$ with a distinguished vertex $v_1$. For $t>1$,
{  we obtain $G_k(t)$} from $G_k(t-1)$ by adding a vertex $v_t$ and a set of $k-1$ edges from
$v_t$ to $G_k(t-1)$ chosen in the following way.
Pick a $k$-clique $C$ of $G_k(t-1)$ uniformly at random ({\em uar}),
and choose a $(k-1)$-dimensional face $F$ of $C$ uar.
Extend $F$ to a $k$-clique $C_t$
by the addition of edges from $v_t$ to the vertices of $F$.
Let the vertex set of $C$ be $\{u_1,...,u_k\}$.
As $\{v_t,u_1,...,u_k\}$ induces no $k$-cliques except $C$ and $C_t$,
the number of $k$-cliques in $G_k(t)$ is $t$.
We are interested in the height of $G_k(t)$ above the distinguished vertex $v_1$, i.e. the maximum graph distance from $v_1$
to any vertex of $G_k(t)$.

\ignore{
Associated with  $G_k(t)$ is   tree $\G_k(t)$ in which each vertex of $\G_k(t)$ corresponds
to a clique of $G_k(t)$. Thus $\G_k(1)$ is the isolated vertex $C_1$, and adding the $k$-clique $C_t$  to $G_k(t-1)$
as described above
 is interpreted as
connecting a leaf vertex $C_t$ to a vertex
$C$ of $\G_k(t-1)$.
The tree $\G_k(t)$ is a useful descriptive device, which we will use
from time to time, but
we emphasize that the height of $G_k(t)$ is not the
same as the height of $\G_k(t)$.
}

A random 2-tree $G_2(t)$ is obtained by joining $v_t$ to a random end point of
a random edge; i.e. by preferential
attachment. Thus the height of random 2-trees is given by the result of
 Pittel \cite{Pit} discussed above.
In the case of a random tree  on $t$ vertices generated by preferential attachment
 Pittel \cite{Pit} established the \whp\ result that the height
$h(t)$ of the tree is asymptotic to $h(t) \sim c \log t$ where $c = 1/(2\g)$ and
$\g$ is the smallest positive solution to $1+ \g + \log \g=0$.
We include this result in our general statement as a special case.
It follows naturally as a special case of our method and serves to check correctness
of the base case.

\begin{theorem}\label{Th1}
For $k \ge 2$ let $h(t ; k)$ be the height of a random $k$-tree on $t$ vertices. Then \whp\
$h(t;k) \sim c \log t$ where $c$ is given as follows:

{\bf Case $k=2$}  \cite{Pit},
$c$ is the solution of
\beq{PITL}
\frac{1}{2c} \exp\set{1+\frac{1}{2c}}=1. \nonumber
\eeq

{\bf Case $k \ge 3$ constant},
$c$ is the solution of
\[
\frac{1}{c} = \sum_{\ell=0}^{k-2} \frac{k}{\ell+ak},
\]
where the value of $a$ is given by
$$
\frac{\G(k) \G(ka)}{\G( ka+k-1)} \; \exp \brac{
\sum_{\ell=0}^{k-2} \frac{ka+k-1}{\ell + ak}}=1,
$$
and $\G(\ell)$ is the gamma function.

{\bf Case $k \rai$}
\[
c \sim  \frac{1}{k \log 2}.
\]
\end{theorem}
In the above theorem, and in  Theorem \ref{Th2} we assume that $k$ is constant or tends slowly to infinity with $t$.
The bound on $k$ used in the proofs is $k=o(\log^{1/3} t)$, but we do not attach any special significance to this value.

The table  below compares asymptotic value of
height (rounded up to the next integer) and results found by experiment
 for  $k$-trees on $t=2^{27}$ vertices.
\begin{center}
\begin{tabular}{|c||c|c|c|c|c|c|c|c|c|c|}
\hline
Value of $k$&2&3&4&5&6&8&10&12&15&20\\
Height: Experimental result&16&10&8&7&5&4&4&3&3&3\\
Height: $\rdup{\log t/(k \log 2)}$&14&9&7&6&5&4&3&3&2&2\\
\hline
\end{tabular}
\end{center}
See Figure \ref{fig1} in the appendix for a
plot of the results obtained as a function of $t$, and Figure \ref{fig2} for the fit to $\rdup{\log t/(k \log 2)}$.


\paragraph{The height of random $k$-Apollonian networks.}
An Apollonian network
is the generalization of an Apollonian triangulation, which can be described as follows.
Initially there is single triangle embedded in the plane.
At the first step this triangle $ABC$ is
divided into three
by insertion of a point $D$ in the interior of the
triangular face and adding lines $DA, DB, DC$.
The triangles
$ABD, ACD, BCD$ replace the original triangle $ABC$
in the embedded triangulation of the plane.
At each subsequent step
some triangular face is subdivided in the same manner.

A random $k$--Apollonian network $A_k(t)$, $t \ge 0$, is obtained as follows.
Start with $A_k(0)$ a $k$-clique
$C_0=K_k$ with vertex set $\{c_1,...,c_k\}$ embedded in $k-1$ dimensions.
For $t>0$,
make $A_k(t)$ from $A_k(t-1)$
 by adding a vertex $v_t$, and edges
chosen as follows.
Pick a $k$-clique $C$ of $A_k(t-1)$ uar. Let the vertex set of $C$ be $U=\{u_1,...,u_k\}$.
Insert $v_t$ in the interior of $C$
and join
$v_t$ to $u_i, i=1,...,k$ by an edge $u_iv_t, \; i=1,...,k$.
This replaces $C$ by $k$ new embedded cliques with vertices $U+v_t-u_i,\; i=1,...,k$.
(We will use the notation $U+v-w$ to mean $(U \cup \{w\}) \setminus \{v\}$.)
The number {  of} embedded $k$-cliques in $A_k(t)$ is $(k-1)t+1$.

\begin{theorem}\label{Th2}
For $k \ge 3$ let $h(t ; k)$ be the height of a random $k$--Apollonian
network on $t$ vertices.
Then \whp\
$h(t;k) \sim c \log t$ where $c$ is given as follows:

{\bf Case $k \ge 3$ constant}
$c$ is the solution of
\beq{dick}
\frac{1}{c} = \sum_{\ell=0}^{k-1} \frac{k-1}{\ell+a(k-1)},
\eeq
where the value of $a$ is given by
\beq{duck}
\frac{k!}{ (a(k-1)) \cdots((a+1)(k-1))}
\exp \brac{ \sum_{\ell=0}^{k-1} \frac{(k-1)(a+1)-1}{\ell+a(k-1)}} =1.
\eeq
{\bf Case $k \rai$}
\[
c \sim  \frac{1}{k \log 2}.
\]
\end{theorem}
Recently, and independently of this work results for random Apollonian networks
were obtained by
Ebrahimzadeh et al \cite{Wor} and Kolossv\'ary \cite{Kol}. The work of \cite{Wor}
adapted the results of Broutin and Devroye
\cite{BroDev} to derive the height of random Apollonian triangulations.
The value of $c=0.8342...$
they obtained is the solution to \eqref{dick}, \eqref{duck} with $k=3$ and corresponds
to a value of $a=2.0683$.
The work of \cite{Kol} uses a different approach based on codes combined with
general techniques for
Markov processes.
In an earlier work, Frieze and Tsourakis \cite{FT} bounded the height of random
Apollonian triangulations
from above
by the height of a random 3-branching, using a result of \cite{BroDev}.

\paragraph{General method}
The technique we describe is simple, bypasses the classic
continuous time branching process results,
(no prior knowledge needed)
and works well for the complicated multi-type branching processes involved in $k$-trees,
Apollonian networks etc.
The main requirement is that the quantities $W_i(t)$ we estimate can be expressed
as recurrences of the form
\[
W_i(t+1) = W_i(t) +\frac{1}{t} \sum_{j \le i} \a_{ij} W_j(t).
\]
By partitioning the steps $t=0,1,2...$ into small intervals, lower and upper bound
approximations for $W_i(t)$ are obtained which can then be expressed via rational
generating functions
from which the coefficients can be extracted.

\section{The height of random $k$-trees}
\subsection{Proof outline of the main theorems} \label{prfout}
In the construction of $G_k(t)$
we add a new vertex  $v_t$ at each step and extend this vertex to a unique $k$-clique.
Starting with a copy of $K_k$ as $G_k(1)$,
the number of of $k$-cliques in $G_k(t)$ is  equal to the number of steps $t$.

 We use the  parameters
$\om=\om(t)$, $s$, and $m$ given by
\beq{assup}
\om=\ooi \log^{2/3}t, \qquad s=\rdup{t^{1/\om}}=o(t), \qquad m = \frac{\log t/s}{\log \brac{1+\frac{1}{\om}}},
\eeq
where the $\ooi$ term  is chosen to make $m$ integer.
The proof in  Section \ref{extract}  assumes $k=o(\sqrt\om)$, but our choice of $\om$ is somewhat arbitrary anyway.
The purpose of the parameter $\om$ is to interpolate the interval $s,s+1,...,t$ at steps $s_j=s (1+1/\om)^j$, $j=0,...,m$.

We will prove Theorem \ref{Th1} in the following way:
\begin{enumerate}

\item Let the height $h(G_k(t))$ be denoted by $h(t)$.
We break our analysis of $h(t)$ into two parts.

Let $C$ be a fixed $k$-clique  added at some step $1,...,s$.
If we consider only those $k$-cliques added at steps $s+1,...,t$, then some (possibly empty)
subset of these form a $k$-tree $G_C(t)$ rooted at $C$, i.e. at the lowest labeled vertex of $C$.
This $k$-tree $G_C(t)$ is
a subgraph of $G_k(t)$.
Note that
if $C, C'$ are distinct $k$-cliques added at steps $1,...,s$
the subgraphs $G_C(t)$ and $G_{C'}(t)$ have no $k$-cliques in common.

Let
the $k$-cliques added at steps $i=1,...,s$ be indexed $C_i: i=1,...,s$.
By the height of the subtree $G_C(t)$, we mean the height of $G_C(t)$ rooted at (the lowest labeled vertex of) $C$.
The main problem is to obtain an asymptotic estimate for the maximum of the heights of the subtrees $G_{C_i}(t), \; i=1,..,s$.
Let
\beq{hst}
 h_s(t) = \max_{C_i, i=1,...,s}\set{ h(G_{C_i}(t))}. 
\eeq
Informally, $h_s(t)$ is the height of $G_k(t)$ if we regard the first $s$ of the $k$-cliques as rooted at level zero.
The fact that $C$ may be an ancestor of $C'$ in
$G_k(s)$ is not relevant to our estimate of $h_s(t)$.

Let $h_0(s)$ be the height $h(s)$ of $G_k(s)$ rooted at $v_1$.
The height $h(t)$ of $G_k(t)$ is bounded by
\beq{host}
h_s(t) \le h(t) \le h_0(s)+ h_s(t)+1.
\eeq

\item In our description of the BFS tree $T_v$ rooted at a
distinguished vertex $v$. Let $N =0,1,...$ denote the levels of $T_v$.
For every $k$-clique $C$ in $G_k(t)$, there is a level $N$, such that the vertices of $C$
lie only in levels $N$ and $N+1$ of $T_v$ (see Section \ref{Recexact}).
We use the  notation
$[N,(\ell, k-\ell)]$ to refer to those cliques with $\ell$ vertices
in level $N$ of the BFS tree $T_v$ and $k-\ell$
vertices in level $N+1$.

\item
Let
$W_{N,\ell}(t)$ be the expected number of $[N,(\ell,k-\ell)]$ configured cliques at step $t$
rooted at any of the first $s$ cliques.
In Section \ref{Recexact} we obtain a recurrence for $W_{N,\ell}(t)$.
In Sections \ref{lowerbd} and \ref{upperbound} we obtain generating functions for
lower and upper bounds $W^L_{N,\ell}(t)\le W_{N,\ell}(t)\le W^U_{N,\ell}(t)$.

\item
Let
$W_N(t)=W_{N,2}(t)$ be the expected number of $[N,(2,k-2)]$ configured cliques at step $t$
rooted at any of the first $s$ cliques.
The height of these cliques above $G_k(s)$ is $N+1$.
As the height depends on $N$ but not $\ell$, the value  $\ell=2$ was chosen for convenience
in the proof.

\item In Section \ref{extract} we see how to extract the coefficients of the generating functions for $W^L_N(t),
 W^U_N(t)$.

\item Let $N= c \log (t/s)$, and let $N'=(1-\e)N$, $N''=(1+\e)N$
for some $\e \ra 0$.
In Sections \ref{seclower} and  \ref{secupper}
we find a value of $N$, and hence $c$ such that
$W^U_{N''}(t) \ra 0$ but $W^L_{N'}(t) \rai$.
Thus w.h.p. $h_s(t) < N''$. The value of $c$ we obtain is the one given in Theorem \ref{Th1}.

\item
In Section \ref{conclb} we prove that the height of a random $k$-tree at step $t$
is at least $N'$ w.h.p.

\item
Let $h_0(s)$ be the height of $G_k(s)$ rooted at $v_1$.
We argued above that
$$h_s(t) \le h(t) \le h_0(s)+ h_s(t).$$
In Lemma \ref{lemup} we prove that $h_0(s)= O( \log s)$ \whp,
thus for $s$ as given in \eqref{assup}, $\log s= (\log t)/\om$.
 We have established that \whp\ $h_s(t) \sim N = c \log (t/s)$,
and thus
\[
 c \log t \le h(t) \le c \log t + O((\log t)/\om).
\]
As we assume that $k^2=O(\om)$, and from our proof $c=O(1/k)$,
it follows
that $G_k(t)$ has height $h(t) \sim c \log t$, \whp
\end{enumerate}

\ignore{
The proof of Theorem \ref{Th2} follows more or less automatically
from  the results obtained for the
proof of Theorem \ref{Th1}. The main difference for random
$k$-Apollonian networks, is  that at step
$t$ the number of vertices $n(t)=t+k$ and the number of embedded cliques is $(k-1)t+1$. }

\subsection{Recurrence for tree height}\label{Recexact}
We will  describe the structure of $G_k(t)$ in terms of
the levels of vertices within each clique relative to the root vertex $v_1$.
The following example using $k=3$ is instructive of our labeling method.
In $G_3(1)$ the initial clique $C_1=K_3$ with vertices $\{v_1,u,u'\}$
We place $v_1$ at level $i=0$ of the BFS tree and $u,u'$ at level $i+1=1$.
The index of the lowest level of $C_1$ is $i=0$
and $C_1$ is oriented $(1,2)$ in that one vertex ($v_1$) is at level $i$
and two vertices ($u,u'$) are at level $i+1$.
We will say that $C_1$ is configured $[0,(1,2)]$. Extending a face of
$C_1$ using vertex $v_2$ gives rise to three possibilities.
If face $v_1u$ or $v_1u'$ is chosen, we obtain another $[0,(1,2)]$ configured clique $C$.
If face $uu'$ is chosen we obtain a $[1,(2,1)]$ configured clique
$\set{u,u',v_2}$ between levels $i=1$ and $i+1=2$.

 We recall the inductive definition of a random $k$-tree.
A random $k$-tree $G_k(t)$, $t \ge k$, is obtained as follows.
Start with $G_k(1)$ a $k$-clique $C_1$ with a distinguished vertex  $v_1$, i.e. with vertices $v_1,x_2,...,x_k$, say.
At subsequent steps $t>1$,
 we obtain $G_k(t)$ from $G_k(t-1)$ by adding a $k$-clique $C_t$ with distinguished vertex $v_t$
and a set of $k-1$ edges from
$v_t$ to $G_k(t-1)$ chosen as follows.
Pick a $k$-clique $C$ of $G_k(t-1)$ uniformly at random.
Let the vertex set of $C$ be $\{u_1,...,u_k\}$.
Choose a $(k-1)$-dimensional face $F$ of $C$ uar.
Suppose, for the purposes of description that the vertices of $F$ are $\{u_1,...,u_{k-1}\}$.
Extend $F$ to a $k$-clique $C_t=\{v_t,u_1,...,u_{k-1}\}$
by the addition of edges from $v_t$ to the vertices of $F$.
As $\{v_t,u_1,...,u_k\}$ induces no $k$-cliques except $C$ and $C_t$,
the number of $k$-cliques in $G_k(t)$ is $t$, and the number of vertices is $t+k-1$.
The precise $k$-cliques which have been added to
form $G_k(t)$ can be found by choosing  the $k$-clique containing the vertex with the highest label $v_t$, and deleting this
vertex recursively.

In general we use the notation {\em clique} to refer to a $k$-clique which has been
added according to our recursive process, and {\em face} to refer to a clique of
dimension $k-1$. We regard $G_k(t)$ as rooted a vertex $v_1$.
We are interested in the height of $G_k(t)$ rooted at $v_1$.
The level sets of the vertices of the breadth first
search tree rooted at $v_1$ form a convenient descriptive device. Inductively the vertices of
each $k$-clique $C$  lie in two adjacent levels $i$ and $i+1$ of this BFS tree.
The notation $[i, (\ell,k-\ell)]$ describes a $k$-clique $C$ with $\ell$ vertices at
level $i$ and $k-\ell$ vertices at
level $i+1$ relative to the BFS tree rooted at $v_1$.
In this case we say $C$ is  $[i, (\ell, k-\ell)]$ configured. In this notation,
the initial clique $C_1$ containing the foot vertex $v_1$
is $[0, (1,k-1)]$ configured.

Given that $C$ is  $[i, (\ell, k-\ell)]$ configured, the number and type of possible extensions of  faces $F$ of $C$
to a new $k$-clique $C'$
are obtained as follows.
An extension of $C=\{u_1,...,u_k\}$ consists of deleting a
vertex $u_j$ (to obtain a face $F$) and then
inserting a vertex $v$
to form $C'=\{ u_1,...,u_{j-1},v,u_{j+1},...,u_k\}$.
If the deleted $u_j$ is chosen
among the $k-\ell$
vertices at level $i+1$
then $C'$ is configured $[i,(\ell,k-\ell)]$. If $u_j$ is chosen among the $\ell$
vertices at level $i$,
then provided $\ell>1$,
$C'$ is configured $[i,(\ell-1,k-\ell+1)$. In the special case $\ell=1$,  $C$ is configured
$[i,(1,k-1)]$, and $u_1$ is the only vertex at level $i$. Deleting $u_1$
results in a clique $C'$ configured $[i+1, (k-1,1)]$ with $k-1$ vertices at level $i+1$ and one vertex at level $i+2$.


Our first step is as follows.
We first obtain bounds for $h_s(t)$ as defined in \eqref{hst}.
Referring to \eqref{host}
we will argue later that, for suitable choice of $s$, we have $h_0(s)=o(h_s(t))$, and hence
$h(t) \sim h_s(t)$.

We will modify the
notation $[i,(\ell,k-\ell)]$ to deal with  our calculation of $h_s(t)$.
Recall that $h_s(t)$ is the height of $G_k(t)$ if we regard the first $s$ of the $k$-cliques as rooted at level zero.
Let $C'$ be a clique that was added at steps $s+1,...,t$. Then $C'$ is a
descendant of some cliques $C$ added at the first $s$ steps,  i.e. some
$C \in \{C_1,...,C_s\}$. In this case we say that $C'$ is
{\em relatively} configured $[i,(\ell,k-\ell)]$ if it is configured $[i,(\ell,k-\ell)]$ in the sub-$k$-tree $G_C(t)$
rooted at $C$. The definition of $G_C(t)$ is given item 1. of Section \ref{prfout}.

Let $W_{i,\ell}(t)$ be the expected number of $[i,(\ell,k-\ell)]$ relatively
configured cliques in $G_k(t)$.
By assumption, at step $s$,  $W_{0,1}(s) =  s$.
We have the following recurrences.
\beq{rec1}
W_{0,1}(s) =  s, \qquad \qquad W_{0,\ell}(t) =  0, \qquad \ell \ge 2, t \ge s.
\eeq
{\bf Case $i=0$: $  [0,(1,k-1)]$ relatively configured cliques.}
\beq{rec2}
W_{0,1}(t+1)=  W_{0,1}(t)+ \frac{k-1}{k} \frac{W_{0,1}(t)}{t}.
\eeq
{\bf Case  $ \ell=k-1,i \ge 1$: $ [i,(k-1,1)]$ relatively configured cliques.}
\beq{rec3}
W_{i,k-1}(t+1)= W_{i,k-1}(t)
+ \frac{1}{k} \frac{W_{i,k-1}(t)}{t} + \frac{1}{k} \frac{W_{i-1,1}(t)}{t}.
\eeq
{\bf Case $ \ell \ne k-1, i \ge 1$: $ [i,(\ell,k-\ell)]$ {relatively configured cliques }.}
\beq{rec4}
W_{i, \ell}(t+1)=  W_{i,\ell}(t)
+ \frac{k-\ell}{k} \frac{W_{i,\ell}(t)}{t} + \frac{\ell+1}{k}
\frac{W_{i,\ell+1}(t)}{t}.
\eeq
The recurrences \eqref{rec2}-\eqref{rec4} can be explained in the following way.
To be specific, consider \eqref{rec4}, and $t>s$. Let $\ul W_{i,\ell}(t)$ be a random variable giving the
number of $[i,(\ell,k-\ell)]$ relatively
configured cliques in $G_k(t)$.
Then
\[
\E (\ul W_{i, \ell}(t+1) \mid \ul W_{i,\ell}(t),\ul W_{i,\ell+1}(t)) =  \ul W_{i,\ell}(t)
+ \frac{k-\ell}{k} \frac{\ul W_{i,\ell}(t)}{t} + \frac{\ell+1}{k}
\frac{\ul W_{i,\ell+1}(t)}{t}.
\]
The term $(k-\ell)/k$ is the probability to pick a face $F$ with $k-\ell-1$ vertices at level $i+1$, from a $[i, (\ell, k-\ell)]$
relatively configured clique. Similarly, the term $(\ell+1)/k$ is the probability to pick a face $F$ with $\ell$ vertices at level $i$
from a $[i, (\ell+1, k-\ell-1)]$ relatively configured clique. Taking expectations again gives \eqref{rec4}.

\subsection{Lower bound for $W_{i,\ell}(t)$: Generating Function}
\label{lowerbd}

As all of the
$W_{i,\ell}(t)$ are monotone non-decreasing in $t$, replacing $t$ by $t'\le t$ in
\eqref{rec1}-\eqref{rec4}
gives a lower bound for the expected number of  $[i,(\ell,k-\ell)]$ relatively
configured cliques at step $t+1$. Recall from \eqref{assup} that $s=\rdup{t^{1/\om}}$.
For $j \ge 0$ we will break the steps $s,s+1,...,t$ into
intervals $I_j = [s_j, s_{j+1}-1]$ where $s_0=s$
and $s_j= \rdup{s \l_0^j}$. Here $\l_0=1+1/\om$
where $\om$ is given by \eqref{assup}.
For fixed $t$ we choose $\l_0$ to ensure
$s_m=s \l_0^m=t$, so that
\beq{mmm}
m = \frac{\log t/s}{\log \l_0}.\nonumber
\eeq
We can assume that $\om$ is chosen so that $m$ is an integer.

We now describe a sub-process which gives lower bounds $W^L\leq W$.
Basically, to do this, for $\t \in I_j$ we replace
$W_{i,\ell}(\t)$ by $W^L_{i,\ell}(s_j)$, so that only vertices which choose cliques
from the  lower bound sub-process added
before $s_j$ count towards the growth of the sub-process.
Thus during $I_j$ the equations corresponding to \eqref{rec1}-\eqref{rec4}
for the sub-process can be replaced by the following.
\begin{flalign}
W^L_{0,1}(s_{j+1}) = & \;W^L_{0,1}(s_{j})+ \frac{k-1}{k} W^L_{0,1}(s_{j})
\sum_{\t = s_j}^{s_{j+1}-1} \frac{1}{\t}\label{W0} \\
W^L_{i,k-1}(s_{j+1})=&\; W^L_{i,k-1}(s_{j}) + \brac{ \frac{1}{k} {W^L_{i,k-1}(s_j)} +
 \frac{1}{k} {W^L_{i-1,1}(s_j)}}\sum_{\t = s_j}^{s_{j+1}-1} \frac{1}{\t}\label{Wi}\\
W^L_{i,\ell}(s_{j+1})=&\; W^L_{i,\ell}(s_{j}) +
\brac{ \frac{k-\ell}{k} {W^L_{i,\ell}(s_j)} +
\frac{\ell+1}{k} {W^L_{i,\ell+1}(s_j)}}
\sum_{\t = s_j}^{s_{j+1}-1} \frac{1}{\t}.\label{Wii}
\end{flalign}
If $f(x)$ is monotone decreasing
\[
f(a+1)+\cdots+f(b) \le \int_a^b f(x) dx \le f(a)+\cdots+f(b-1).
\]
Thus
\[
\sum_{\t = s_j}^{s_{j+1}-1} \frac{1}{\t} - \frac{1}{s_j} \le \int_{s_j}^{s_{j+1}}
\frac{dx}{x} \le \sum_{\t = s_j}^{s_{j+1}-1} \frac{1}{\t}.
\]
As $s_j= \rdup{s\l_0^j}$ it follows that
\begin{flalign}\label{do1}
\sum_{\t = s_j}^{s_{j+1}-1} \frac{1}{\t}
=& \frac{\th_1}{s_j} + \log \frac{\rdup{s\l_0^{j+1}}}{\rdup{s\l_0^j}}\\
=& \log \l_0(1+\d_j),\label{do2}
\end{flalign}
where $0 \le \th_1 \le 1$ and $|\d_j| \le 2/s_j$ provided $s \rai$.

Substitute \eqref{do1}--\eqref{do2} for the summation in \eqref{W0}--\eqref{Wii}.
Let $\d'= \max_j |\d_j|$, thus $\d'=o(1/\om)$ (see \eqref{assup}). Let
\beq{LAMD1}
\l_1=\l_0 (1-\d')= \l_0(1-o(1/\om)).
\eeq
Replace $\l_0$ with $\l_1= \l_0(1-\d')$ to obtain a uniform lower bound on the recurrences for all $j$,
and re-scale by dividing by $s$
to obtain simplified recurrences
$W^L_{i,\ell}(j)\leq W_{i,\ell}(s_j)/s$.
We obtain
\begin{align}
W^L_{0,1}(0)= &1,\label{low1}\\
W^L_{0,\ell}(0)= &\; 0 \qquad \ell \ge 2,\label{low2}\\	
W^L_{0,1}(j+1) = &\; W^L_{0,1}(j) \brac{1 + \frac{k-1}{k} \log \l_1},\label{low3}\\
W^L_{i,k-1}(j+1) = &\; W^L_{i,k-1}(j) \brac{1+\frac{1}{k} \log \l_1}+ W^L_{i-1,1}(j)
\frac{1}{k} \log \l_1,
\quad i\geq 1,\label{low4}\\
W^L_{i,\ell}(j+1)= &\; W^L_{i,\ell}(j) \brac{1+ \frac{k-\ell}{k} \log \l_1} +
W^L_{i,\ell+1}(j) \frac{\ell+1}{k} \log \l_1,\quad{  i\geq 1,\ell\neq k-1}.\label{low5}
\end{align}
Let $G^L_{i,\ell}(z)$ be the generating function
for $W^L_{i,\ell}(j), \; j\ge 0$, and let $\g_{\ell}= 1+ ((k-\ell)/k ) \log \l_1$.
It follows from \eqref{low1}, \eqref{low3} that
\[
G^L_{0,1}(z) =  \frac{1}{1-\g_1 z}.
\]
From \eqref{low2}, \eqref{low4}, \eqref{low5}, we obtain
\beq{geffy0}
G^L_{i,k-1}(z)= \g_{k-1} z\;  G^L_{i,k-1}(z) + \brac{\frac{1}{k} \log \l_1 } z \;
G^L_{i-1,1}(z), \nonumber
\eeq
\beq{geffy}
G^L_{i,\ell}(z)= \g_{\ell} z\;  G^L_{i,\ell}(z) + \brac{\frac{\ell+1}{k} \log \l_1 }
z \;G^L_{i,\ell+1}(z),
\quad{  i\geq 1,\ell\neq k-1}. \nonumber
\eeq
Thus
\begin{flalign}
G^L_{i,k-1}(z)=& \frac{1}{k}\frac{z  \log \l_1 }{1-\g_{k-1}z}\; G^L_{i-1,1}(z),
\quad i\geq 1,\label{GEF1}\\
G^L_{i,\ell}(z) = & \frac{\ell+1}{k} \frac{z \log \l_1 }{1-\g_{\ell}z} \; G^L_{i,\ell+1}(z)
\quad{  i\geq 1,\ell\neq k-1}.\nonumber
\end{flalign}
It follows inductively that
\beq{geffy2}
G^L_{i,1}(z) = \bfrac{ z^{k-1}k!(\log \l_1)^{k-1}}{k^k (1-\g_1z)\cdots(1-\g_{k-1}z)}^i
\frac{1}{1-\g_1z},
\eeq
and for $\ell = 2,...,k-2$
\beq{geffy3}
G^L_{i,\ell}(z) = \frac{1}{k} \; \prod_{j=\ell}^{k-1} \frac{j+1}{1-\g_j z}
\bfrac{z \log \l_1}{k}^{k-\ell}G^L_{i-1,1}(z).
\eeq

\subsection{Upper bound for $W_{i,\ell}(t)$: Generating Function}\label{upperbound}

For simplicity of notation, put
$\a_\ell=(k-\ell)/k$
and $\b_\ell =(\ell+1)/k$. Then iterating the main variable backwards in recurrences \eqref{rec2} --
\eqref{rec4}, and recalling that $W_{i,\ell}(t)$ is non-decreasing in $t$ gives
\begin{align*}
W_{0,1}(t+\s)&=  W_{0,1}(t)\prod_{j=0}^{\s-1}\brac{1+ \frac{\a_1}{t+j}}\\
W_{i,k-1}(t+\s)&= W_{i,k-1}(t)\prod_{j=0}^{\s-1}\brac{1+ \frac{\a_{k-1}}{t+j}}
+ \a_{k-1} \sum_{j=0}^{\s-1} \frac{W_{i-1,1}(t+j)}{t+j}
\prod_{i=j+1}^{\s-1}\brac{1+\frac{\a_{k-1}}{t+i}}\\
&\leq W_{i,k-1}(t)\prod_{j=0}^{\s-1}\brac{1+ \frac{\a_{k-1}}{t+j}}
+ \a_{k-1} \sum_{j=0}^{\s-1} \frac{W_{i-1,1}(t+\s)}{t+j}
\prod_{i=j+1}^{\s-1}\brac{1+\frac{\a_{k-1}}{t+i}}\\
W_{i,\ell}(t+\s)&= W_{i,\ell}(t) \prod_{j=0}^{\s-1} \brac{1+\frac{\a_\ell}{t+j}} +
\b_\ell \sum_{j=0}^{\s-1} \frac{W_{i,\ell+1}(t+j)}{t+j}
\prod_{i=j+1}^{\s-1}\brac{1+\frac{\a_\ell}{t+i}}\\
&\leq W_{i,\ell}(t) \prod_{j=0}^{\s-1} \brac{1+\frac{\a_\ell}{t+j}} +\b_\ell
\sum_{j=0}^{\s-1}\frac{W_{i,\ell+1}(t+\s)}{t+j}
\prod_{i=j+1}^{\s-1}\brac{1+\frac{\a_\ell}{t+i}}
\end{align*}
Let $t=s_j$, let $ t+\s= s_{j+1}$ and let $W_{i,\ell}(j)=W_{i,\ell}(s_j)/s$ for all $i,j,\ell$.
Thus
\begin{align*}
W_{0,1}(j+1)&=W_{0,1}(j)\prod_{t=s_j}^{s_{j+1}-1} \brac{1+\frac{\a_1}{t}}\\
&\leq  \sbrac{1+O\bfrac{1}{s}}\l_0^{\a_1}W_{0,1}(j)\\
&\leq  \brac{1+\frac{k-1}{k}\log\l'}W_{0,1}(j).
\end{align*}
For $\l_0=1+1/\om$ the value of $\l'=1+1/\om +O(1/\om^2)$. To see this, for $a<1$
 the function $f(x)=x^a-(1+a \log x)$
has a unique minimum at $x=1$, with $f(1)=f'(1)=0$, so the Taylor expansion of $f(1+h)=O(h^2)$.

Similarly
\begin{align*}
W_{i,k-1}(j+1)&\leq W_{i,k-1}(j)\prod_{t=s_j}^{s_{j+1}-1}
\brac{1+\frac{\a_{k-1}}{t}} +\a_{k-1}W_{i-1,1}(j+1) \sum_{t=s_j}^{s_{j+1}-1}
\frac{1}{t}\prod_{\t=t+1}^{s_{j+1}-1}\brac{1+\frac{\a_{k-1}}{\t}}\\
&\leq  \sbrac{1+O\bfrac{1}{s}} \brac{ W_{i,k-1}(j)
\l_0^{\a_{k-1}} + W_{i-1,1}(j+1)(\l_0^{\a_{k-1}}-1)}\\
&\leq W_{i,k-1}(j) \brac{1+\frac{1}{k} \log \l'}+ W_{i-1,1}(j+1)\frac{1}{k}\log \l' ,
\end{align*}
and
\begin{align*}
W_{i,\ell}(j+1) &\leq  W_{i,\ell}(j)
\prod_{t=s_j}^{s_{j+1}-1} \brac{1+\frac{\a_\ell}{t}} +\b_\ell
W_{i,\ell+1}(j+1) \sum_{t=s_j}^{s_{j+1}-1}
\frac{1}{t}
\prod_{\t=t+1}^{s_{j+1}-1}\brac{1+\frac{\a_\ell}{\t}}\\
&\le \sbrac{1+O\bfrac{1}{s}} \brac{ W(j) \l_0^{\a_\ell} +\frac{\b_\ell}{\a_\ell} W_{i,\ell+1}(j+1)
(\l_0^{\a_\ell}-1)}\\
&\le W_{i,\ell}(j) (1+\a_\ell \log \l')+\b_\ell W_{i,\ell+1}(j+1)\log \l' .
\end{align*}
We thus
obtain the following recurrences for an upper bound
$W^U_{i,\ell}(j)\geq W_{i,\ell}(s_j)/s$.
\begin{align*}
W^U_{0,1}(0)= &1,\\
W^U_{0,\ell}(0)= &\; 0 \qquad \ell \ge 2,\\	
W_{0,1}^U(j+1)&=\brac{1+\frac{k-1}{k}\log\l'}W_{0,1}^U(j).\\
W_{i,k-1}^U(j+1)&=W_{i,k-1}^U(j) \brac{1+\frac{1}{k} \log \l'}+
W_{i-1,1}^U(j+1)\frac{1}{k}\log \l', \quad i\geq 1,\\
W^U_{i,\ell}({j+1})& = W^U_{i,\ell}(j)\brac{1+\frac{k-\ell}{k}\log \l'}+
\frac{\ell+1}{k} W^U_{i,\ell+1}({j+1}){  \log\l'},\quad{  i\geq 1,\ell\neq k-1}.
\end{align*}
Let $G^U_{i,\ell}(z)$ be the generating function
for $W^U_{i,\ell}(j), \; j\ge 0$, and let $\g'_{\ell}= 1+ ((k-\ell)/k ) \log \l'$.
It follows that
\[
G^U_{0,1}(z) =  \frac{1}{1-\g'_1 z},
\]
and generally, we obtain
\beq{geffy0x}
G^U_{i,k-1}(z)= \g'_{k-1} z\;  G^U_{i,k-1}(z) + \brac{\frac{1}{k} \log \l' }
\;G^U_{i-1,1}(z), \nonumber
\eeq
\beq{geffyx}
G^U_{i,\ell}(z)= \g'_{\ell} z\;  G^U_{i,\ell}(z) + \brac{\frac{\ell+1}{k} \log \l' }
 \;G^U_{i,\ell+1}(z),
\quad{  i\geq 1,\ell\neq k-1}. \nonumber
\eeq
Thus
\begin{flalign}
G^U_{i,k-1}(z)=& \frac{1}{k}\frac{\log \l' }{1-\g'_{k-1}z}\; G^U_{i-1,1}(z),
\quad i\geq 1,\label{GEF1X}\\
G^U_{i,\ell}(z) = & \frac{\ell+1}{k} \frac{\log \l' }{1-\g'_{\ell}z} \; G^U_{i,\ell+1}(z)
\quad{  i\geq 1,\ell\neq k-1}.\nonumber
\end{flalign}
Choosing $\l'_1=\max{\l'}$, it follows inductively that
\beq{geffy2x}
G^U_{i,1}(z) = \bfrac{k!(\log \l'_1)^{k-1}}{k^k (1-\g'_1z)\cdots(1-\g'_{k-1}z)}^i
\frac{1}{1-\g'_1z},
\eeq
and for $\ell = 2,...,k-2$
\beq{geffy3x}
G^U_{i,\ell}(z) = \frac{1}{k}\; \prod_{j=\ell}^{k-1} \frac{j+1}{1-\g'_j z}
\bfrac{\log \l'_1}{k}^{k-\ell}G^U_{i-1,1}(z).
\eeq
The expressions \eqref{GEF1X}, \eqref{geffy2x}, \eqref{geffy3x}
differ from \eqref{GEF1}, \eqref{geffy2}, \eqref{geffy3}
in that $\l'$ replaces $\l$ and multiplicative powers of $z$ are suppressed.
It will be seen in Section \ref{extract}
that the $W^L$ and $W^U$ are sufficiently close to obtain tight bounds on the
the expected occupancy of each level.

\section{Random $k$-trees: Asymptotic expression for maximum height}\label{extract}

We now show how to
extract the coefficients of our generating functions.
For reasons of symmetry of the generating function it is easier for us to focus on
$W^X_{N,2}(t)$ for $X=L,U$ and suitable $N \rai$.
Choosing $\ell=2$  will suffice. The height of the rooted $k$-tree depends on the index $N$ of the maximum  level set, and not  on $\ell=1,...,k$. As the case $k=2$ is already known from \cite{Pit}, we assume $k \ge 3$.

\subsection{Extraction of coefficients for a lower bound on the expected size of level sets}
\label{seclower}

We first discuss the case for
$G^L(z)=G^L_{N,2}(z)$, where from \eqref{geffy2} and \eqref{geffy3}
\beq{GN2}
G^L(z) = \frac{k}{2 z \log \l_1}\bfrac{ z^{k-1}k! (\log \l_1)^{k-1}}{k^k(1-\g_1z)
\cdots(1-\g_{k-1}z)}^N. \nonumber
\eeq
Using  $[z^m]G(z)$ to denote the coefficient of $z^m$ in the formal expansion of $G(z)$,
let $w^L(m)=W^L_{N,2}(m)=[z^m]G^L(z)$.
To extract the coefficient, let $M= m-N(k-1)-1$, so that
\beq{remember}
[z^m]G^L(z) = \frac{k}{2  \log \l_1}
\brac{(\log \l_1)^{k-1}k!/k^k }^N [z^M] (f(z))^N,
\eeq
where
\beq{fofx}
f(z) = \frac{1}{(1-\g_1z)\cdots(1-\g_{k-1}z)}.
\eeq
We want the smallest $N$ such that $[z^m]G^L(z) \ra 0$, i.e. $W^L_{N,2}(t) \ra 0$.
We will write $N$ as  $N =c \log t$.
We  prove in the next paragraph that $c=\Th(1/k)$.
Given these restrictions on $c$, we have $N/m \ra 0$ and $N=N(t) \rai$.

By inspection, as $\g_\ell >1$, the coefficients of
$(f(z))^N$ are at least as large as the coefficients of $1/(1-z)^{N(k-1)}$. Recall that $m=(\log t/s )/\log \l_0$. Thus
\begin{align*}
[z^m] G^L(z) & \ge \frac{k}{2  \log \l_1}
\brac{(\log \l_1)^{k-1}k!/k^k }^N { N(k-1)-1+M \choose M}\\
&\ge \frac{k}{2  \log \l_1}
\brac{(\log \l_1)^{k-1}k!/k^k }^N \bfrac{m-2}{ N(k-1)-1}^{N(k-1)-1}\\
&= \Th\bfrac{k^2 N}{\log t}  \brac{\frac{\sooi\log (t/s)} {N(k-1)}\bfrac{k!}{k^k}^{1/(k-1)}          }^{N(k-1)}.
\end{align*}
As $N \ge 1$ and assuming $k \ge 3$, the value of $[z^m] G^L(z)$ tends to infinity with $t$ for any $Nk \le a\log(t/s)$
for some constant $a$. In particular, if $k \rai$ then we can choose $a=1/2e$.

For an upper bound, we have
\beq{maxel}
 \max_{\ell} \g_{\ell}= \max_{\ell} (1+ (k-\ell)/k \; \log \l_1) = \g_1= 1+\frac{k-1}{k} \log \l_1.
\eeq
Using $M= m-N(k-1)-1$, and $\log \l_1=\log \l_0(1+o(1/\om))$ from \eqref{LAMD1},
\beq{MonN}
\frac{N}{M}= (1+O(1/\om))c\log \l_0= c' \log \l_1.
\eeq
It follows that
\[
\g_1^M \le e^{M \frac{k-1}{k} \log \l_1} \le \brac{e^{\sooi/ck}}^{N(k-1)}.
\]
Thus
\begin{align*}
[z^m] G^L(z) & \le \frac{k}{2  \log \l_1}
\brac{(\log \l_1)^{k-1}k!/k^k }^N { N(k-1)-1+M \choose M}\g_1^M\\
&\ge \frac{k}{2  \log \l_1}
\brac{(\log \l_1)^{k-1}k!/k^k }^N \bfrac{me}{ N(k-1)}^{N(k-1)-1}\g_1^M\\
&= \Th\bfrac{k^2 N}{\log t}  \brac{\frac{\sooi\log (t/s)} {N(k-1)}e^{1/ck}        }^{N(k-1)}.
\end{align*}
Thus $[z^m] G^L(z)$ tends to zero with $t$ for any $k \ge 3$ and  $c \ge 3/k$.

We next describe a general technique (based on \cite{Kin}) to obtain an asymptotic expression for $[z^M](f(z))^N$
in terms of
an implicitly defined parameter $\wh x$.
The method can be broken into six steps.
\begin{enumerate}[{\bf M1}]
\item Write
\[
[z^M] (f(z))^N= \frac{f(x)^N}{x^M} [z^M] \bfrac{f(zx)}{f(x)}^N.
\]
\item Let $Y(x)$ be a random variable with  probability generating function
 $\E z^Y= f(zx)/f(x)$.
 By inspection of the generating function, (see \eqref{fofx}) the random variable $Y$
has positive probabilities on the non-negative integers.
 Let $Y_1,...,Y_N$ be i.i.d. as $Y$.
 \[
 [z^M] \bfrac{f(zx)}{f(x)}^N = [z^M] \E (z^{Y_1+ \cdots+ Y_N})=
\Pr (Y_1+\cdots+ Y_N=M).
 \]
\item
Obtain the moments $\mu(x), \s^2(x)$ of $Y(x)$ from
\begin{align*}
&\mu(x) = \E Y =  \left. \frac{d}{dz }\E z^Y \right|_{z=1}= x \frac{f'(x)}{f(x)}\\
&{  \s^2(x)-\m(x)+\m(x)^2=\E Y(Y-1)=  \left. \frac{d^2}{dz^2}\E z^Y
\right|_{z=1}= x^2\frac{f''(x)}{f(x)}}
\end{align*}
\item Choose $Y$ so that $\mu(x)=\E Y = M/N$.
Solve $\mu(x)= M/N$ for $x$.
\item
We have chosen $\E (Y_1+\cdots+ Y_N)=M$. Provided $\s^2(x) $ is bounded, and
as the random variable $Y$
has lattice width $h=1$, by the Local Limit Theorem (see e.g. \cite{Dur} or \cite{Gned})
\beq{LLT}
\Pr(Y_1+\cdots+ Y_N=M) = (1+O(1/N)) \frac{1}{\sqrt{2\pi \s^2N}}. \nonumber
 \eeq
\item
From \bM{1}, \bM{2} and \bM{5},
\[
[z^M] f(z)^N=(1+O(1/N)) \frac{1}{\sqrt{2\pi \s^2N}} \frac{f(x)^N}{x^M}.
\]
The value of $x$ is obtained from the condition that $\mu(x) = M/N$  in {\bf M4},
 and the value of $\s^2(x)< \infty$ from {\bf M3}.
\end{enumerate}
We apply this method  to $f(z)$ from \eqref{fofx}. For step {\bf M3} we find
\begin{eqnarray}
\mu(x)& = &
 \sum_{\ell=1}^{k-1} \frac{\g_{\ell} x}{1-\g_{\ell}x} \label{muu}\\
\s^2(x) &=& \mu(x) +
 \sum_{\ell=1}^{k-1} \frac{(\g_{\ell} x)^2}{(1-\g_{\ell }x)^2} \label{sigma}.
 \end{eqnarray}
 Considering {\bf M4}, from \eqref{MonN} we can relate $\mu(x)$ to $c'$ by
 \beq{muxc}
 \mu(x) = \frac{M}{N}= \frac{1}{c' \log \l_1}.
 \eeq
Recall that $\g_{\ell}= 1+ (k-\ell)/k \; \log \l_1$, and that $\l_1=(1+1/\om+ o(1/\om^2))$.
From \eqref{maxel} $ \max_{\ell} \g_{\ell}=\g_1$.
 In order to find the value of $\wh x$ from \eqref{muu} note that
the smallest singularity of \eqref{muu} is at $x=1/\g_1=1-O(1/\om)$.
 Intuitively, as $M/N \rai$ it must be that $\wh x \ra 1/\g_1$
from below.
From \eqref{muu} it follows that if $\mu(x)>0$ then  $ x >0$.
The function $g(x) =  \sum_{\ell=1}^{k-1} (\g_{\ell} x)/(1-\g_{\ell}x)$
is monotone increasing in $x$
from $g(0)=0$.
Thus the solution $x>0$  to $g(x)=\mu(x)$ is unique.

 Based on this, for some $a=a(k)$, to be determined, let
 \beq{whx}
 \wh x = \frac{1-a \log \l_1}{\g_1}=\frac{1-a \log \l_1}{1+((k-1)/k) \log \l_1}
=1- O\bfrac{a+1}{\om}
 \eeq
From
$$\frac{1}{c'\log\l_1}=\m(\wh x)\leq \frac{k\g_1{\wh x}}{1-\g_1{\wh x}}
\leq\frac{k}{a\log\l_1},$$
it follows that $a>0$. Also as we assumed  $c'=\Th(1/k)$, we have that $a=O(1)$.
From these observations and \eqref{whx}
we see that
\beq{xo1}
{\wh x}=1-O(1/\om).
\eeq
The next step is to prove $a=\Th(1)$ so that \eqref{500} (see below) is bounded, and establish
the relationship between $a$ and $c$ given in \eqref{asymp}.
From \eqref{muu} we have
\begin{flalign}
\mu(\wh x) = & \frac{1-a \log \l_1}{\log \l_1}
\brac{\sum_{\ell = 1}^{k-1}
\frac{1+\frac{k-\ell}{k} \log \l_1}{\frac{\ell-1}{k}+a+\frac{a(k-\ell)}{k }\log \l_1}}\label{first1}\\
=&\brac{1+O\bfrac{k}{\om}}\frac{k}{\log \l_1}\; \sum_{\ell=0}^{k-2}
\frac{1}{\ell+ak}. \label{last}
\end{flalign}
The right hand side of \eqref{first1} is monotone decreasing in $a$ from (at least) $kH_{k-1}/\log\l_1$ when $a=0$.
Thus if  $\mu(\wh x)=M/N=(1+O(k/\om))/c \log \l_0$, and $c=\Th(1/k)$ there is a unique $a$ satisfying the relationship $c=c(a,k)$.
We see from \eqref{last} that
\beq{asymp}
\frac{1}{kc}=\brac{1+O\bfrac{k}{\om}}\sum_{\ell=0}^{k-2} \frac{1}{\ell+ak}.
\eeq
Note that
\beq{bint}
\int_x^{x+j+1} \frac{dy}{y}  \le \frac{1}{x}+\cdots+\frac{1}{(x+j)}
\le \frac{1}{x}+ \int_{x}^{x+j} \frac{dy}{y}.
\eeq
Putting  $x=ak$ we see that
\beq{sumakl}
\log \frac{k(a+1)-1}{ka} \le \sum_{\ell=0}^{k-2}
\frac{1}{\ell+ak} \le \frac{1}{ka}+\log \frac{k(a+1)-2}{ka-1}.
\eeq
This implies that
\[
\frac{1}{kc} \ge \ooi \log \frac{k(a+1)-1}{ka},
\]
so that
\[
a \ge \frac{1-1/k}{e^{\sooi/ck}-1}.
\]
As we assumed $c=\Th(1/k)$, it follows for $k \ge 2$ that $a$ is bounded below by a positive constant,
and thus $a=\Th(1)$.
This bound on $a$ combined with \eqref{sumakl} implies that
\beq{500}
\sum_{\ell=0}^{k-2} \frac{1}{\ell+ak}= \log \frac{k(a+1)-1}{ka} +\frac{\z_1}{k}
=\log\frac{a+1}{a}+\frac{\z_2}{k}
\eeq
where $|\z_1|,|\z_2|=O(1)$.

Thus crudely, for some $B=\Th(1)$
\beq{Bmu}
\mu(\wh x)= B\frac{k}{\log\l_1}. \nonumber
\eeq

Armed with this, our next task is to approximate $\s^2(\wh x)$, as given in \eqref{sigma}.
Writing $\s^2(x) = \mu(x) + \f(x)$ and substituting \eqref{whx} we find,
for some $B'=\Th(1)$ that
\[
{  \f(\wh x)} = B' \frac{k}{\log^2 \l_1}.
\]
Thus
\beq{sigx}
\s^2(\wh x) = B\frac{k}{\log\l_1}+B' \frac{k}{\log^2 \l_1}. \nonumber
\eeq
Proceeding to step {\bf M6} and using $m/N=1/(c \log \l_0)$, $M=m-N(k-1)-1$ and
\eqref{remember} we have that
\begin{flalign}
W^L_{N,2}(m)=
[z^m] G^L(z) = &\frac {k \wh x}{2 \log \l_1} \frac{1+O(1/N)}{\sqrt{2 \pi \s^2N}}
\left[\frac{k!}{k^k} (\log \l_1)^{k-1} \frac{f(\wh x)}{\wh x^{1/c\log \l_0}}
\wh x^{k-1} \right]^N
 \nonumber\\
=& \Th\bfrac{k^{1/2}}{N^{1/2}}\;\; [\F(k,a)]^N. \label{Phika}
\end{flalign}

The final step is to put $\F(k,a)$ into a more tractable form by removing
the parameter $c=c(a,k)$. Our aim is to prove
\beq{bunny}
\F(k,a)=\brac{1+O\bfrac{k^2}{\om}}
\frac{\G(k) \G(ka)}{\G( ka+k-1)} \; \exp \brac{
\sum_{\ell=0}^{k-2} \frac{ka+k-1}{\ell + ak}}.
\eeq
This can be done as follows:
\begin{enumerate}[{\bf F1.}]
\item
From the definition of $\wh x$ in \eqref{whx}
\beq{sothat}
\wh x^{-1/c \log \l_0} = \exp \brac{ \frac{1}{kc} \brac{ka+ {k-1}}
\brac{1+O({1}/{\om})}}.
\eeq
\item
It follows from
\eqref{asymp} that
\[
\frac{1}{kc}=\brac{1+O\bfrac{k}{\om}}\sum_{\ell=0}^{k-2} \frac{1}{\ell+ak}.
\]
As $c=\Th(1/k)$ the right hand side sums to a constant, so that \eqref{sothat} can be written as
\[
\wh x^{-1/c \log \l_0} = \brac{1+O\bfrac{k^2}{\om}}
 \exp \brac{\sum_{\ell=0}^{k-2} \frac{ka+k-1}{\ell + ak}}.
\]
\item From the definition of $\wh x$
\beq{xerr}
\wh x^{k-1} = 1+ O\bfrac{k}{\om}.\nonumber
\eeq
\item From \eqref{fofx} and the definition of $\wh x$
\begin{eqnarray*}
\frac{(\log \l_1)^{k-1}}{k^{k-1}} f(\wh x) &=&
\brac{1+ \frac{k-1}{k} \log \l_1}^{k-1} \prod_{{  \ell=1}}^{k-1}
\frac{1}{{  \ell-1+ak+a(k-\ell)\log \l_1}}\\
&=&  \frac{1+O(k/\om)}{ (ka)(ka+1)\cdots(ka+k-2)}.
\end{eqnarray*}
\end{enumerate}
Putting F1 to F4 together gives us \eqref{bunny}.

\subsection{Extraction of coefficients for an upper bound on the expected size of level sets}
\label{secupper}
We now consider $w^U(m)=W^U_{N,2}(m)=[z^m]G^U(z)$ where $G^U(z)=G^U_{N,2}(z)$.
Observe first that
\beq{similar}
G^U_{N,2}(z)=\frac{k}{2\log \l_1'}\brac{\frac{k!}{k^k}\frac{(\log \l_1')^{k-1}}{(1-\g_1'z)\cdots(1-\g_{k-1}'z)}}^N
\nonumber
\eeq
Then with $m=M$ as in Section \ref{seclower} we have from \eqref{remember} that
and $\l'=\l_1',\g_\ell'$ replacing $\l_1,\g_\ell$ in $f$,
\begin{align}
w^U(m)=[z^m]G^U(z)& = \frac{k}{2  \log \l'}
\brac{(\log \l')^{k-1}k!/k^k }^N [z^{M}] f(z)^N\nonumber\\
&\leq \frac{k}{2  \log \l'}\brac{(\log \l')^{k-1}k!/k^k }^N \frac{f({\wh x})^N}
{{\wh x}^{M'}}\nonumber\\
&=\Theta\bfrac{k^{1/2}}{N^{1/2}}\brac{\brac{1+O\bfrac{k}{\om}}\F(k,a)}^N,\label{wLwU}
\end{align}
where $\F(k,a)$ is given in \eqref{Phika}.

\subsection{Asymptotic value of maximum height}

For a given height $N=c \log(t/s)$, we get lower and
upper bounds for $W_{N,2}(t)$, the expected number of $[N,(2,k-2)]$ relatively configured $k$-cliques
from   \eqref{Phika} and \eqref{wLwU}, which in turn depend on $\F(k,a)$.
Provided  $w^L(m) \sim w^U(m)$,  and we can prove concentration of the level set sizes around  these bounds,  the maximum height $h_s(t)$
can be obtained
 from the value  $a$ which makes $\F(k,a)=1$ in \eqref{bunny}.
By expanding $\F(k,a)$ around this value of $a$ we prove that the value $w^U(m) \ra 0$ for larger values of $c$, whereas
$w^L(m) \rai$ for smaller values of $c$.

Once we find the
value of $a$ such that
  $\F(k,a)=1$, we can obtain $c(a)$ via \eqref{asymp}.
Our analysis of behavior around $\F(k,a)=1$ depends on $k$.
Basically there are three cases. $k=2$, $k\ge 3$ constant, and $k \rai$.

Setting aside the details for now,  for $k \ge 3$ constant, the implicit relationships \eqref{bunny} and \eqref{asymp}
 are the content of Theorem \ref{Th1}.
In the case that  $k \rai$, the value of $a$ solving $\F(k,a)=1$
can be obtained explicitly as  $a=1+o(1)$. This solution allows us to
obtain an explicit asymptotic of $c\sim 1/k \log 2$ from \eqref{asymp}.

\subsubsection{Case $k \rai$}
We assume $k=o(\sqrt{\om})$ and use the asymptotic expansion of \eqref{bunny}.
As $\G(y)= e^{O(1/y)} {  \sqrt{2 \pi }} \;\; y^{y-1/2} e^{-y}$
 $\F(k,a)$ can be written as
\begin{align}
\F(k,a)= &
\brac{1+O\bfrac{1}{k}+O\bfrac{k^2}{\om}} \;\;
\bfrac{(a+1)^3 2 \pi k}{a}^{1/2}\nonumber
\\
& \times \brac{
\frac{a^a}{(a+1)^{a+1}}
 \exp \brac{(a+1-1/k)\sum_{\ell=0}^{k-2} \frac{1}{ak+\ell}}}^k. \label{rhsval}
\end{align}

It is easiest to expand directly about $a=1$. The value of $\F(k,1)$ is
\beq{Fk1}
\F(k,1) \sim (16 \pi k)^{1/2} \frac{1}{4^k} (2-\b)^{2k-1}, \nonumber
\eeq
where $0\le \b\le 1/k$.
This value $\b$ is deduced as follows.  Putting $a=1$, we can bound the sum in \eqref{rhsval} by
\beq{bint}
\log \frac{2k-1}{k}= \int_k^{2k-1} \frac{dy}{y}  \le \frac{1}{k}+\cdots+\frac{1}{2k-2}
\le \int_{k-1}^{2k-2} \frac{dy}{y} = \log \frac{2k-2}{k-1}. \nonumber
\eeq
We see that for some $0 \le \b \le 1/k$,
\beq{logb}
\sum_{\ell=0}^{k-2}\frac{1}{k+\ell} = \log (2-\b).
\eeq
Denote the final bracketed term on the RHS of \eqref{rhsval} above by
$\Psi(k,a)^k$.
Note that
\beq{aeq1}
\Psi(k,1) = \frac{1}{4} (2-\b)^{2-1/k}.
\eeq
The expansion of $\F(k,\a)$ in $\a=a(1+\e)$ can thus be obtained by
expanding $\Psi(k,\a)$ about $a=1$. We write $\Psi(k,a(1+\e))= F_1 e^{F_2}$, and
note the following simplifications.
\begin{align}
F_1= &\frac{ (a(1+\e))^{a(1+\e)}}{(a(1+\e)+1)^{a(1+\e)+1}}=
\frac{a^a}{(a+1)^{a+1}} \bfrac{a^a}{(a+1)^a}^{\e}
\frac{(1+\e)^{a(1+\e)}}{(1+a\e/(a+1))^{a+1+\e a}} \label{Fo}\\
=& \frac{a^a}{(a+1)^{a+1}} \bfrac{a^a}{(a+1)^a}^{\e}
e^{O(\e^2)}.
\label{F1}
\end{align}
The second line comes from an expansion of the last
term on the right hand side of \eqref{Fo},
using $(1+x)=\exp (\log (1+x))$,
in which the first order terms disappear.
\begin{align}
F_2=&(a(1+\e)+1-1/k) {  \sum_{\ell=0}^{k-2}} \frac{1}{ak(1+\e)+\ell} \nonumber\\
 = &(a+1-1/k) {  \sum_{\ell=0}^{k-2}} \frac{1}{ak+\ell}
+\e \brac{\sum_{\ell=0}^{k-2} \frac{a}{ak+\ell}
- a(k(a+1)-1) {  \sum_{\ell=0}^{k-2}}
\frac{1}{(ak+\ell)^2}} +O(\e^2).
\label{second}
\end{align}
Thus,
\beq{fefi}
\Psi(k, a(1+\e))=
 \Psi(k,a)\;\; \brac{ e^{O(\e)}\bfrac{a}{a+1}^a \exp \brac{- a\sum_{\ell=0}^{{  k-2}}
 \frac{k-\ell-1}{(ak+\ell)^2}}}^{\e},
\eeq
which, for $\e>0$ decreases faster than $(a/(a+1))^{a\e}$.

Applying \eqref{bint} with $b=2$ and $a=1$ to the second term in \eqref{second} gives
\[
(2k-1) \sum_{\ell=0}^{k-3} \frac{1}{(k+\ell)^2 }= 1-\th,
\]
where $\th =O(1/k)$.
As a result, from \eqref{logb},  \eqref{aeq1}, \eqref{F1}, \eqref{second} and \eqref{fefi}
\beq{rhsval0}
\Psi(k,1+\e)= \Psi(k,1) \brac{e^{O(\e)} \brac{1-\frac{\b}{2}}e^{-1+\th}}^{\e }=
\frac{1}{4}(2-\b)^{2-1/k} \brac{ e^{O(\e)} \brac{1-\frac{\b}{2}}e^{-1+\th}}^{\e}.
\eeq
The coefficients $w^L(m), \; w^U(m)$ we wish to evaluate at $a=1+\e$ are given in \eqref{Phika} and \eqref{wLwU},
respectively by
\begin{align}
w^L(m)=W^L_{N,2}(m)&= \Th\bfrac{k^{1/2}}{N^{1/2}} [\F(k,a)]^N\label{wNm}\\
w^U(m)=W^U_{N,2}(m)&\leq
\Th(k)\left[\brac{1+O\bfrac{k}{\om}}\F(k,a)\right]^N\label{wNm0}.
\end{align}
 From \eqref{rhsval} and \eqref{rhsval0} and $\th = O(1/k)$.
\[
\F(k,1+\e)=\brac{1+O\bfrac{1}{k}+O\bfrac{k^2}{\om}+O(\e)+O(k\e^2)} \;\;
\brac{4 \pi k}^{1/2}
\brac{\brac{1-\frac{\b}{2}}^{2+\e}
\; e^{-\e}}^k.
\]
The the $O(1/k)$ is from $1/(1-\b/2)$, the
 $O(\e)$ is from $e^{\th \e k}$ and the $O(k \e^2)$ from $e^{k O(\e^2)}$.
These come from raising the expression in \eqref{rhsval0} to the power $k$.

\paragraph{Upper bound on height.}
Choose $|\e|= A(\log k)/k$ for some constant $A>0$,
then  as $0 \le \b\le 1/k$, for some $0 \le \xi \le 1$,
\[
\F(k,1+\e)= (2+o(1)) \sqrt{\pi} e^{-\xi} \sqrt{k}e^{-\e k}=k^{1/2-A+o_k(1)}\text{ as }
k\to\infty.
\]
Suppose first that $A>1/2$, say $A=1$. Then for large enough $k$,
$$W_{N,2}(t)\leq sw^U(m)\leq s\Theta(k)k^{-N/3}\to 0.$$
Here we use $k=o(\sqrt{\om})$ (see below\eqref{assup}) to ensure convergence to zero.

We show in Lemma \ref{lemup}  that \whp\ the
height of $G_k(\t)$ is bounded by $O(\log\t)$. It follows that \whp\ the height
\beq{upheight}
h(t;k)\leq O(\log s)+N+1.
\eeq
Indeed, we have shown that w.h.p. there are no $[N,(2,k-2)]$
relatively configured cliques
and the clique generation process means that in this case there will be no
$[N+1,(\ell,k-\ell)]$ relatively configured cliques.

From \eqref{muxc} we have $\m(\wh x)\sim \frac{1}{c\log\l_1}$. From
\eqref{last} and \eqref{500} we have
$\m(\wh x)\sim \frac{k}{\log\l_1}\cdot\log\frac{a+1}{a}$.
It follows that $c\sim \frac{1}{k\log 2}$.
From \eqref{upheight} we have  w.h.p. that
\beq{64a}
h(t;k)\leq O\bfrac{\log t}{\om}+c\log(t/s)+1\sim \frac{\log t}{k\log 2}.
\eeq
This proves the upper bound in Theorem \ref{Th1} for the case where $k\to\infty$.

\paragraph{Lower bound on height.}
Now consider the lower bound. Putting $A<0$ we get from
\eqref{wNm} that
$$W_{N,2}(m)\geq w^L(m)\geq \Th\bfrac{k^{1/2}}{N^{1/2}}k^{-AN}\to\infty.$$
We show in Section \ref{conclb} that this is good enough
to prove that $h(t;k) \ge (1-o(1)) c \log t$ w.h.p. This establishes
a lower bound asymptotic to \eqref{64a}. Thus as asserted by Theorem \ref{Th1}
\[
 h(t;k) \sim \frac{\log t}{k\log 2}.
 \]
\subsubsection{Case $k$ constant}
The case $k=2$ can be resolved by our methods, but it is  proved
in \cite{Pit} and the paper is already long enough, we omit this case.
For $k \ge 3$, the statement of Theorem \ref{Th1} follows from
 \eqref{bunny}, and the following details.
\[
\frac{\G(k) \G(ka)}{\G( ka+k-1)}= \frac{(k-1)!}{(ka+k-2)(ka+k-3)\cdots(ka)}.
\]
Let $a$ be the unique positive solution to $\F(k,a)=1$. Let $\a=a(1+\e)$, then
\begin{align}
(k\a + k-2)\cdots(k\a) = & (ka+k-2)\cdots(ka)
\prod_{\ell= 0}^{k-2}\brac{1+ \frac{\e ka}{ka+\ell}}
\nonumber\\
=&(ka+k-2)\cdots(ka)\exp\brac{ \e k a \sum_{\ell=0}^{k-2}
\frac{1}{ka+\ell}+ O(\e^2 k)}.\label{last1}
\end{align}
Using \eqref{last1} to deal with the exponential term in the definition
of $\F(k,a)$ in  \eqref{bunny}, we see that
\[
\F(k, \a) = \F(k,a) \exp\brac{
- a(k(a+1)-1)  \sum_{\ell=0}^{k-2}
\frac{1}{(ak+\ell)^2}
}^{\e k}\times e^{O(\e^2k)}.
\]
We now see from \eqref{wNm0} that if $\e>0$ and $\e N\to \infty$ then
$w^U(m)\to 0$ and so from \eqref{64a} we see that w.h.p.
$h(t;k)\leq (1+o(1))c\log t$ where the value of $c$ is given by \eqref{asymp}.
This verifies the upper bound in Theorem \ref{Th1} for this case.

When $\e<0$ and $-\e N\to \infty$ we see from
\eqref{wNm} that $w^L(m)\to\infty$. In
Section \ref{conclb} we show that $w^L(m)\to\infty$ suffices to prove  with high probability
that $h(t;k)\geq (1-o(1))c\log t$. This verifies the lower bound in Theorem \ref{Th1} for this case.

\section{Concentration of occupancy of level sets around expected value $W_{i,\ell}(t)$}
\label{conclb}

The coefficient $W^L_{N,2}(t)$ is the expected
value of a random variable $\ul W$
corresponding to a subprocess of $G_k(t)$. Recall that
$h_s(t)$ is our estimate of the expected height of $G_k(t)$ above $G_k(s)$.
If we can prove concentration of $\ul W$ from below
for $H=(1-\e)h_s(t)$, then the height of $G_k(t)$
is at least  $H$ \whp\
To do this we follow
the method of Devroye  \cite{Dev1a}, which we translate
into our discrete step context. This method
couples the growth of the level sets with a suitably defined Galton-Watson process.
We first explain our approach.
Because we observe the process at a given step $t$
 the total number of vertices added is fixed, and the proof requires an additional twist.

It is convenient to consider coupling our discrete process with a continuous time
process. To do this, we replace the step parameter $t$ of the previous sections by $n$
and reserve  variables such as $t, T, \t$ for times in the continuous process.

Our basic view of the discrete process starting from the clique set $S$ of
$G_k(s)$ is as a set of bins $C_1,...,C_i,...,C_s$. At step $s$ each bin $C_i$
contains a single ball $v_i$, {  corresponding to a single clique}.
Suppose that at step $n \ge s$
bin $C_i$ contains $\nu_i$ balls. At the next step, step $n+1$,
the probability that ball $v_{n+1}$ goes into bin $C_i$ is
$\nu_i/n$. Given the occupancy $\nu_i$ of $C_i$ we can subsequently
 construct a branching  $\cT(\nu_i)$ rooted at clique $v_i$
as a $k$-tree process of length $\nu_i$.

As mentioned above, we wish to use the  method in Devroye  \cite{Dev1a} to
prove concentration of the lower bound.
The main problem for us, is that
 the occupancies of the bin system $C_S=(C_1,...,C_i,...,C_s)$
in the discrete process are not
independent. Let $\nu_i$ be the occupancy of $C_i$ then $\nu_1+\nu_2+\cdots+\nu_s=n$.
Using a continuous time  device we construct  independent
sub-processes which occur in $C_S$ \whp

To avoid confusion between the continuous time and
discrete processes in the subsequent discussion
we adopt the following notation. The discrete process at step $n$
is a system of balls in bins. The continuous time process at time $t$,
is a system of particles in cells. For the continuous time system consisting of
particles $C=\{b_1,b_2,\ldots,\}$, each particle $b\in C$
divides independently into $b,b'$ with waiting time $X_b$ a random variable
with (negative) exponential distribution
rate parameter $\r=1$. If the continuous time
system is observed at time $T$ and contains $n$
particles (i.e. we have $C=\{b_1,...,b_n\}$) then:
\begin{itemize}
\item[(i)] The probability $p_j$ that $b_j$ is the next particle
to divide is $p_j=1/n$.
\item[(ii)] The waiting time from $T$ to the  division event of particle $b_j$
is independent exponential with rate parameter $\r=1$.
\item[(iii)] The rate parameter for the next division in
the entire system of $n$ particles is $\r_n=n$.
\end{itemize}
These results follow from the memoryless properties of the exponential distribution.

A pure birth process of this type is known as a
Yule process, see Feller \cite{Fe}.
Given an initial population of $\th$ particles in a cell $C$ at time
$t=0$, the population $\Pi_\th(\t)$ of $C$ at time $\t$
has  distribution $P_n(\t)=\Pr(\Pi_\th({\t})=n)$ given by
\beq{YP}
P_n({\t})= {n-1 \choose n-\th}e^{-\th {\t}}\brac{1-e^{-{\t}}}^{n-\th}.
\eeq
This is the probability of $k=n-\th$ failures
and $r=\th$ successes in $n$ Bernoulii trials, where there is
a success on the $n$th trial. The probability of success is $p=e^{-{\t}}$.
The expected number of failures $k$ before the $r$-success is $r(1-p)/p$. Thus
\beq{ent}
\E \Pi_\th({\t}) = \th+\frac{\th(1-p)}{p}=\th e^{\t}.
\eeq

In our case the cell $C$ can be regarded either as a single cell $C_S$ with
$\th=s$ at $t=0$, or as
$s$ sub-cells with $\th=1$ at $t=0$; the latter
corresponding to the balls in bins system of the discrete process. By choosing a time
$\t_n=\log(n/\th),$
from \eqref{ent}
 the expected size of the population
is $n$.
We use this relationship to switch between the discrete and
the continuous time processes.
If we observe a given cell $C$ at time ${\t}$
and $C$ has occupancy $N$ then the rooted branching $\cT({\t})$ is
identical with $\cT(N)$
in the discrete process.
It we start at time 0 with a single cell $C$ with occupancy $\th=1$, and stop at
time ${\t}$
with occupancy $\Pi({\t})=\Pi_1(\t)$, we can restart identically distributed processes
$C_1,...,C_{\Pi({\t})}$ stopping at $2{\t}$, and so on. We now
fix our attention on a given cell
$C$ with $\th=1$ at $t=0$.

In the discrete process, choose $\l=e^{1/\om}$
so that $s=n^{\log \l}=n^{1/\om}$. Here we will assume that \eqref{assup}
holds with $t$ replaced by $n$.
Let $s_j=s\l^{jL}, j=0,1,...$ where {  $L=(1/2)\log (n/s)$}.
Now fix $\t=s\l^L$. For a given bin $C_i$, after $\t$ steps
 the expected occupancy is $\nu=\t/s$, where
\beq{nuval}
\nu=\t/s=\l^L=e^{L/\om}=\bfrac{n}{s}^{1/2\om}=s^{(1/2)(1-1/\om)}\sim \sqrt{s}.
\eeq
In the corresponding continuous time process, let
\[
T_j= jL \log\l + \log s {  =\log s_j}
\]
so that
\[
T_{j+1}-T_j = L \log \l = \log (\t/s)=T, \text{ say.}
\]
Intuitively $T_j$ is the equivalent of $s_j$, and $T$ is the equivalent of $\t$.
For a cell starting with $\th=1$ particles, from \eqref{ent}, \eqref{nuval}
\[
\E \Pi(T) = e^T=\t/s \sim \sqrt{s}.
\]
Because of the memoryless property we restart the Yule processes at $T_j$, $j=0,1,...$,
assigning $i=1$ particles per cell. Starting at $T_j$ each
cell grows independently up to $T_{j+1}$, etc.

A cell $C$ is {\bf good}, if after time $T$ has elapsed,
\begin{itemize}
\item[(i)] The occupancy $\Pi(T) \ge \nu$,
\item[(ii)] The branching constructed on the first $\nu=\t/s$ particles in the cell has height
at least $h=c(1-\e) \log \n$ where $\e=o(1)$.
\end{itemize}
If $C$ is good, let $\wt W_h$ be the occupancy of
level $h$ in this sub-process based on the branchings of the first $\nu$ particles, otherwise let $\wt W_h=0$.
In this way we define a Galton-Watson process
with population sizes $X_j,\;j\ge 0$ as follows.
$X_0=1$, $X_1=\wt W_h$ and in general $X_{j+1}$
is the progeny of the surviving particles at level $j$.
Thus if $X_j=\xi$ then $X_{j+1}= X_{j,1}+\cdots+X_{j,\xi}$
where $X_{j,\ell}, \ell=1,...,\xi$ are
independently distributed as $X_1$.

\[
\E \wt W_h \ge \Pr( \Pi(T) \ge \n) \times \wh W
\]
where
$\wh W=W_{h,2}^L(\n)$ is a lower bound
on the expected number of cliques (balls) at height $h$ at time $\n$
defined in Section \ref{lowerbd}. There is the caveat that
$s$ is replaced by $s'=s^{o(1)}$, chosen so that $s^{o(1)}\to\infty$ with $s$.
We run the discrete process to generate the first $\nu$ balls in the box
(particles in the cell), starting the branching from a base set of $s'$ balls
as in Section \ref{lowerbd}.

 In \eqref{YP}, let $\th=1$, replace $\t$ with $T$ and $n$ with $\nu$. Then
\begin{eqnarray}
\Pr( \Pi(T) \ge \n) &=& \sum_{N \ge \nu} P_N(T)\nonumber\\
&=& \sum_{N \ge \nu} \frac{s}{\t}\brac{1-\frac{s}{\t}}^{N-1}\nonumber\\
&=& \brac{1-\frac{s}{\t}}^{\t/s-1}\ge \frac{1}{2e}.\nonumber 
\end{eqnarray}
If we choose $c(a)$ so that the RHS of \eqref{wNm} tends to infinity then we have
\[
\E \wt W_h \ge \wh W/ 2e> 1.
\]
In the associated Galton-Watson process we have $\mu=\E X_1= \E \wt W_h >1$.
For a Galton-Watson process with mean $\mu >1$,
the probability of ultimate survival is $1-q$ where
$q<1$ is the smallest solution of $q=F(q)$.
Here $F(x)$ is the probability generating function of $X_1$.
Let $M=\max X_1$.
As the entire branching was constructed from $\nu$ particles
it must be that $M\leq\n$. We  use a technique from Devroye \cite{Dev1a} to upper bound $q$ by
\[
q \le 1- \frac{\mu}{M}.
\]
Thus
\beq{q}
q\leq 1-\frac{1}{\n}.
\eeq

Let $\s=\log(n(1-\d)/s)$ for $\d=o(1)$ .
Observing the  population $\Pi_S(\s)$ of the complete $s$-cell Yule process  $C_S$
at {  time} $\s$, from \eqref{ent} we have
\[
\E \Pi_S(\s) = s e^{\s} = {  n(1-\d)}.
\]
Let $N=\Pi_S(\s)$ be the population of the complete process at time $\s$, and  let
$A$ be the event that $N \le n$.
We will establish in Lemma \ref{lemmy} below that $\Pr(\overline A)= o(1)$.

Let $B$ be the event that the height $H$
of $\cT(N)$ satisfies
\[
H \ge h\frac{ \s}{T} = c(1-\e) \log \nu \; \frac{\log (n(1-\d)/s) }{\log\n}=
 c(1-\e') \log n/s
\]
where $\e'=o(1)$. Consider the complementary event
$\overline B$  that none of the $s$ independent
Galton-Watson branching processes  survives past
generation $\rdown{\s/T}$.
From \eqref{nuval} $\nu \sim \sqrt{s}$, and using \eqref{q} we have
\[
\Pr(\overline B) \le q^s \le e^{-(1-o(1)) \sqrt{s}}=o(1).
\]
If the event $A$ occurs, then  $N \le n$
and the corresponding tree $\cT(N)$ is a subtree of $\cT(n)$. Thus
\[
\Pr(\text{height of }\cT(n) \ge (1-\e)c \log n/s)
\ge 1-\Pr(\overline A)-\Pr(\overline B) = 1-o(1).
\]
Finally observe that $\log(n/s)\sim \log n$ and
this completes the proof  for the lower bound on height $h_s(n)$.

\begin{lemma}\label{lemmy}
Let $\s=\log (n(1-\d)/s)$. Provided $\d \ge \sqrt{(K \log n)/s}$, and $s =o(\sqrt n)$ we have
\[
\Pr(\overline A) = \sum_{N \ge n+1} P_N(\s) = O(n^{-(K-3)/2}).
\]
\end{lemma}
\proofstart
From \eqref{YP}, with $\th=s$, and $\t=\s$ and $n=N$, we have
\beq{pns}
P_N(\s)= {N-1 \choose N-s}e^{-s \s}\brac{1-e^{-\s}}^{N-s}.
\eeq
Thus for $N \ge n+1$
\begin{align*}
\frac{P_{N+1}}{P_N} & = \frac{N}{N-s+1} \brac{1-\frac{s}{n(1-\d)}}\\
& \le 1+ s \brac{\frac{1}{N-s} - \frac{1}{n(1-\d)}}\\
& \le 1- \frac{s \d}{2n(1-\d)}.
\end{align*}
Thus,
\[
\sum_{N\ge n+1} P_N = O \bfrac{n}{s \d} P_{n+1}.
\]
However, from \eqref{pns}
\begin{align*}
P_{n+1}& = O\brac{\frac{e}{1-\d}e^{-1/(1-\d)}}^s\\
&= O\brac{e^{-s\d^2/2}}.
\end{align*}
Thus
\begin{align*}
\sum_{N \ge n+1} P_N & = O \brac{\frac{n}{s \d}e^{-s\d^2/2}}\\
&= O\brac{n^{-(K-3)/2}}.
\end{align*}
\proofend

\subsection{Upper bound on height}
\begin{lemma}\label{lemup}
The height $h(t)$ of a random $k$-tree $G_k(t)$ is $O(\log t)$ \whp
\end{lemma}
\proofstart
A crude calculation suffices to establish a \whp\ upper bound  of $O(\log t)$.
Consider a shortest path $P(t)=v_t,u_1,...,u_i,v_1$ back from $v_t$ to the
root vertex $v_1$ in $G_k(t)$. As half of the cliques
$C=K_{k}$ in $G_k(t)$ were added by time $t/2$,
\[
\Pr(v_t \mbox{ chooses a clique } C \mbox{ in }G_k(t/2))\geq  \frac{1}{2}.
\]
Thus the expected distance to the root must be (at least) halved by the edge $v_t u_1$.
Whatever the label $s$ of $u_1=v_s$, this halving occurs independently at the next step.
This must terminate \whp\ after $c \log t$ steps, for some suitably large constant $c$,
as we now prove.

Suppose $v_t$ is at height $h$. The number of halving steps in $P(t)$ is at most $\log_2 t$.
Either the $h$ trials  resulted in less than $h/3$
halving steps, or, if not, then $h\leq 2h/3+\log_2 t$.
If $v_t$ is at height $h=c \log t$ then the second case is a contradiction for $c>3$.
The probability of at most $h/3$ halving steps is at most $ e^{-h/144}$.
Putting $h= 300\log t$,
the probability that some vertex has height $h$ is at most $te^{-h/144} \le 1/t \ra 0$.
Thus \whp\ the height of $G_k(t)$ is $O(\log t)$.
\proofend

\section{Random  Apollonian networks}\label{RAN}
We are interested in the height of $A_k(t)$ rooted at vertex $c_1$.
Once again the height of $A_k(t)$ is the maximum distance of a vertex from the root.
The first problem is to describe the structure of $A_k(t)$ relative to this BFS tree.
The following example using $k=3$ is instructive of our labeling method. In $A_3(0)$,
let the initial clique $C_0$
be a triangle with  vertex set $\{a,b,c\}$. Assume
vertex  $a$ is at level 0 of the BFS tree and $b,c$ at level 1.  The index of the
lowest level of $C_0$ is $i=0$
and $C_0$ is oriented $(1,2)$
giving a $[0, (1,2)]$ configured triangle.
Insertion of a vertex $v$ in the interior of $abc$ replaces this triangle by three
new triangles $abv, acv, bcv$.
Triangles $abv, acv$ are configured $[0, (1,2)]$ and $bcv$ configured  $[1, (3,0)]$ in
that all three vertices of this
triangle lie in level $i=1$ of the BFS tree. Once a clique has been subdivided, it is no longer
considered as part of the Apollonian network. In the above example triangle $abc$ is no longer
available for subdivision. To distinguish this case,
we call the cliques available for subdivision {\em embedded}.

In general, suppose clique
$C=K_k$ is configured
$[i, ( \ell , k-\ell)]$ with vertex set $\{ u_1,...,u_{\ell}, v_{\ell+1},...v_k\}$.
If $\ell = 2,...,k$ then inserting a vertex $w$ in the interior of $C$, removes $C$ and
 produces $\ell$ cliques of type $[i, (\ell-1,k-\ell+1)]$
and $k-\ell$ cliques of type $[i, (\ell,k-\ell)]$. If $\ell = 1$, then
insertion of a vertex
inside a clique of type $[i,(1,k-1)]$
forms one clique of type $[i+1, (k, 0)]$ and $k-1$ cliques of type $[i, (1,k-1)]$.

At each step $k$ embedded cliques are created and one is discarded, as it has been subdivided.
Thus, as proved above Theorem \ref{Th2} the number of embedded cliques in $A_k(t)$ is $(k-1)t+1$.
This leads to the following recurrences for the expected number
$W_{i, \ell}(t)$  of $[i, (\ell, k-\ell)]$
configured cliques at step $t$.
\[
W_{0,1}(0)=1, \qquad W_{i,\ell}(0) = 0 \text{ otherwise}.
\]
\[
W_{0,1}(t+1)= W_{0,1}(t) + \frac{k-2}{(k-1)t+1} W_{0,1}(t).
\]
\beq{winkl}
W_{i,k}(t+1)= W_{i,k}(t) -\frac{1}{(k-1)t+1}W_{i,k}(t) + \frac{1}{(k-1)t+1}W_{i-1,1}(t).
\eeq
For $1 \le \ell \le k-1$,
\[
W_{i,\ell}(t+1)= W_{i,\ell}(t){  +} \frac{k-\ell-1}{(k-1)t+1}W_{i,\ell}(t) +
\frac{\ell+1}{(k-1)t+1}W_{i,\ell+1}(t).
\]

\subsection{ Solution of recurrences}
The system of recurrences for $W_{i,\ell}(t)$ and their solution is
very similar to the case for $k$-trees.
We give an outline description only, pointing out where
differences arise.
The main difference is that \eqref{winkl} contains a negative term. However,
as \eqref{winkl} can be rewritten as
\[
W_{i,k}(t+1)= W_{i,k}(t)\brac{1 -\frac{1}{(k-1)t+1}} + \frac{1}{(k-1)t+1}W_{i-1,1}(t),
\]
the lower bound substitution of $W_{i,\ell}(s_j)$ for $W_{i,\ell}(t)$ is still valid.
We obtain (e.g.) the following system of lower bound recurrences, in place of
\eqref{low1} -- \eqref{low5}.
\begin{align*}
W^L_{0,1}(0) & = 1\\
W^L_{0,1}(j+1) & = W^L_{0,1}(j) \brac{1+ \frac{k-2}{k-1} \log \l}\\
W^L_{i,k}(j+1)& = W^L_{i,k}(j) \brac{1-\frac{1}{k-1}
\log \l} +W^L_{i-1,1}(j)\frac{1}{k-1} \log \l\\
W^L_{i,\ell}(j+1)& = W^L_{i,\ell} (j) \brac{1+
\frac{k-\ell-1}{k-1} \log \l} + W^L_{i,\ell+1}(j) \frac{\ell+1}{k-1} \log \l,
\qquad 1 \le \ell\le k-1.
\end{align*}
For $\ell=1,...,k$ let
\[
\g_{\ell} = 1 + \frac{k-\ell-1}{k-1} \log \l.
\]
The lower bound generating functions satisfy
\begin{align*}
G_{0,1}(z)& = \frac{1}{1-\g_1z}\\
G_{i,k}(z)& = {  \frac{1}{k-1} \frac{z\log \l}{1-\g_k z}G_{i-1,1}(z)}\\
G_{i,\ell}(z)& = {  \frac{\ell+1}{k-1} \frac{z\log \l}{1-\g_{\ell}z}G_{i,\ell+1}(z)},
\end{align*}
leading to
\[
G_{i,1}(z) =
 \brac{\frac{z^k k! (\log \l)^k}{(k-1)^k} \frac{1}{(1-\g_1z)\cdots(1-\g_k z)}}^i \frac{1}{1-\g_1 z}.
\]
Work with $G_{N,2}(z)= G_N(z)$ as before, where
\begin{flalign*}
G_{N,2}(z)= &\frac{k-1}{2} \frac{1-\g_1 z}{z \log \l} G_{N,1}(z)\\
=&
\frac{k-1}{2} \frac{1-\g_1 z}{z \log \l}
 \brac{\frac{z^k k! (\log \l)^k}{(k-1)^k} \frac{1}{(1-\g_1z)\cdots(1-\g_k z)}}^N \frac{1}{1-\g_1 z}.
\end{flalign*}
Making  substitutions  $M=m-kN+1$,
$\wh x = (1-a \log \l)/\g_1$ and so forth leads to the following expression for
$\F(k,a)$ (to be compared with \eqref{bunny}),
\[
\F(k,a) = \frac{k! \G(a(k-1))}{ \G((a+1)(k-1)+1)}
 \exp\brac{((k-1)a+k-2)\sum_{\ell=0}^{k-1} \frac{1}{\ell+a(k-1)}},
\]
subject to the asymptotic identity (to be compared with \eqref{asymp}),
\[
\frac{1}{c(k-1)}= \sum_{\ell=0}^{k-1} \frac{1}{\ell+a(k-1)}.
\]
In the case that $k \rai$ we can expand about $a=1$ in a manner identical to $k$-trees to obtain
the asymptotic height
\[
h(t;k) \sim  \frac{\log t }{k \log 2}.
\]
The case $k \ge 3$ constant, is similar.

The concentration of the upper bound follows easily and the concentration of the lower bound from
methods similar to Section \ref{conclb}.  The main difference is that, in the continuous time model,
on division a particle $b$ is replaced by $k \ge 3$ progeny, as opposed to two progeny, which was the case for $k$-trees.

Let $N$ be the number of (multiple) births, in a Yule process in which each particle has $k$ children and then dies,
i.e. there is
an overall increase in population of $k-1$ per birth.
If the original population size is $s$ at time $t=0$, then the population size after $N$-th birth is $\beta_N=(k-1)N+s$.
The probability of $N$ births having occurred by time $t$ is given by
\[
p_N(t)= \prod_{i=1}^N \frac{(k-1)(i-1)+s}{(k-1)i} \; \; \times e^{-st} (1-e^{-(k-1)t})^N.
\]

\newpage

\begin{figure}
\begin{center}

\includegraphics[scale=0.4]{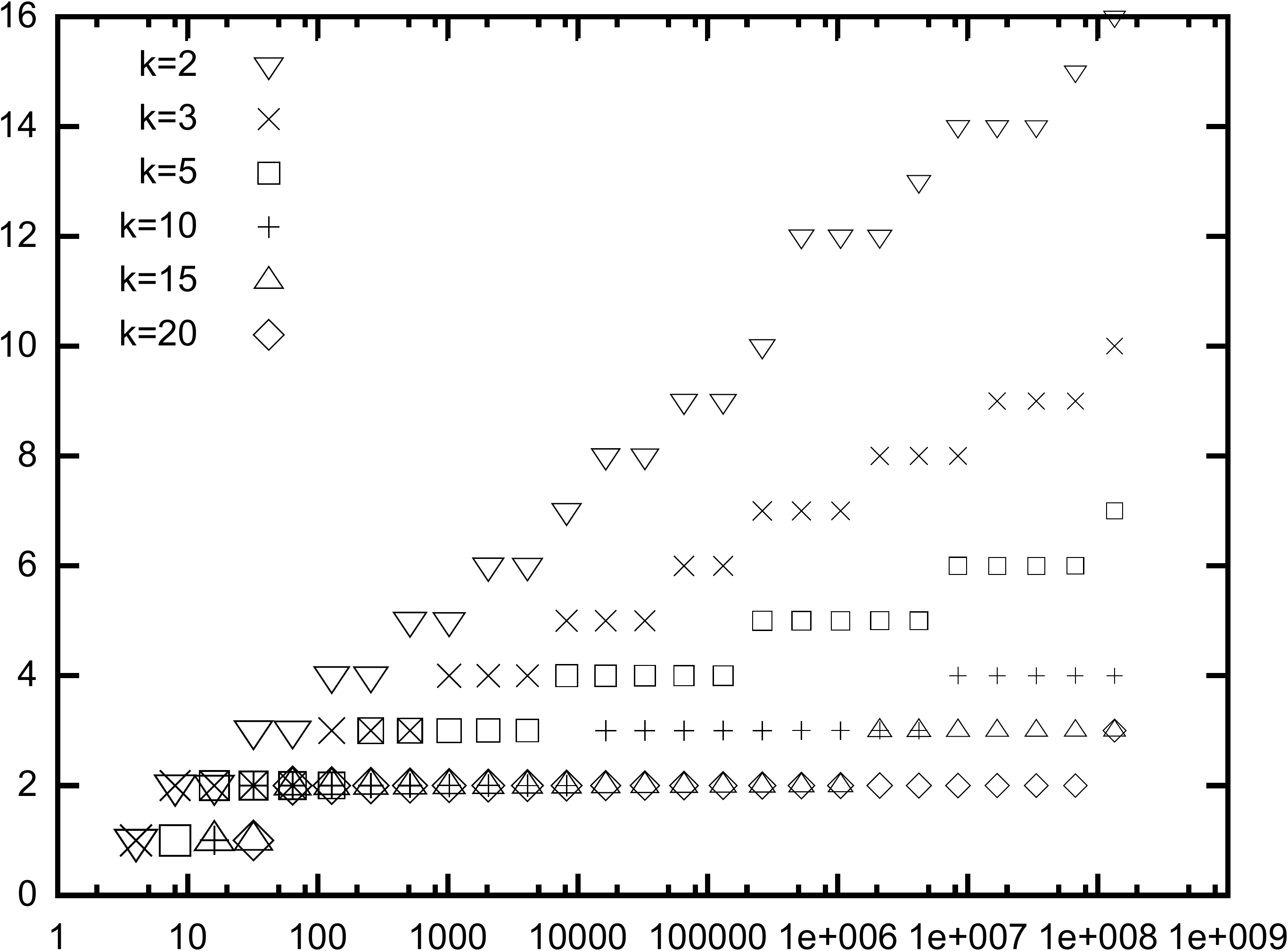}
\end{center}
\caption{Experimental results for the height of random $k$-trees for
$k=2,3,5,10,15,20$}
\label{fig1}
\end{figure}

\begin{figure}
\begin{center}

\includegraphics[scale=0.4]{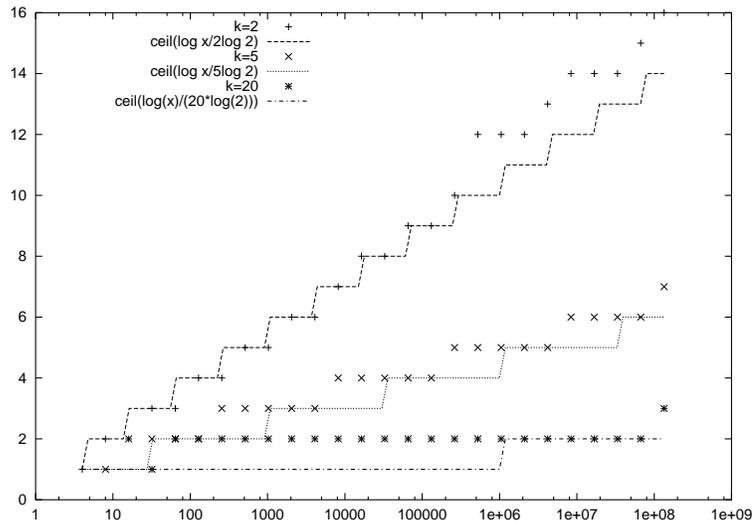}
\end{center}
\caption{Experimental results for $k$-tree height fitted to $\rdup{\log(t)/(k \log 2)}$ for
$k=2,5,20$}
\label{fig2}
\end{figure}

\end{document}